\documentclass[8pt]{article}
\usepackage{graphicx}
\usepackage{caption}
\usepackage{amssymb}
\usepackage{amsmath}
\usepackage{float,color}
\usepackage{amsthm}
\usepackage{subfigure}
\usepackage{mathrsfs}
\newtheorem{theorem}{Theorem}[section]

\newtheorem{lemma}{Lemma}[section]
\newtheorem{proposition}{Proposition}[section]

\newtheorem{remark}{Remark}[section]
\newtheorem{corollary}{Corollary}[section]

\numberwithin{equation}{section}

\DeclareMathOperator{\Fix}{Fix}

\usepackage{physics}

\newcommand{\paren}[1]{\left ( #1 \right )}
\newcommand{\brac}[1]{\left \{ #1 \right \}}

\makeatletter
\newcommand*{\defeq}{\mathrel{\rlap{%
			\raisebox{0.3ex}{$\m@th\cdot$}}%
		\raisebox{-0.3ex}{$\m@th\cdot$}}%
	=}
\makeatother

\usepackage{hyperref}
\hypersetup{colorlinks=true,linkcolor=black,urlcolor=blue,citecolor = black}

\graphicspath{ {images/} }

\renewcommand{\mathrm}{\text}

\usepackage{lineno}
\newcommand{\parcial}[3][ ]{\frac{\partial^{#1} #2}{\partial #3^{#1}}}

\begin{document}

	\title{\bf Mosquito population suppression models with seasonality  and d-concave equations}
	\author{{Rodrigo Fernández-Martínez}\footnote{Departamento de Matemáticas,
			University of Oviedo, Spain.}\;\; and {Alfonso Ruiz-Herrera}\footnote{Departamento de Matemáticas,
			University of Oviedo, Spain (\href{mailto:ruizalfonso@uniovi.es}{ruizalfonso@uniovi.es}).}
	}
	\maketitle
	\begin{abstract}
			The sterile insect technique has emerged recently as a biologically secure and effective tool for suppressing wild mosquito pests. To improve the performance of this strategy, understanding the interaction between wild and sterile mosquitoes is critical. Although the common models for this biological problem are scalar equations, they are remarkably resistant to the mathematical analysis. In a series of papers, Dueñas, Nuñez, and Obaya have developed a powerful approach to describe the dynamical behavior of scalar equations with d-concave nonlinearities, a property typically related to the sign of the third derivative. In this paper, we show that, for periodic equations coming from population dynamics, this condition is typically associated with the positive sign of the third derivative of the inverse of the Poincaré map. This remark allows us to simplify some arguments in the periodic case and obtain a deep geometrical understanding of the global bifurcation patterns. Consequently, the dynamical behavior of the models is analyzed in terms of simple and testable conditions. Our methodology allows us to describe precisely the dynamical behavior of the common mosquito population suppression models, even incorporating seasonality. This paper generalizes and improves many recent results in the literature.
	\end{abstract}
	
	% ----------------------------------------------------------------
	\section{Introduction}
	Mosquito-borne diseases are provoked by bacteria, parasites, or viruses transmitted by mosquitoes. These diseases include chikungunya, dengue, malaria, yellow fever, and Zika \cite{WHO}. Nowadays, mosquitoes have been ranked as one of the world's most dangerous animals to humans. Particularly, dengue fever, transmitted mainly by {\it Aedes aegypti} and {\it Aedes albopictus}, is a public health problem in tropical and subtropical regions where more than 390 million people are annually at risk of dengue infection \cite{Bhatt2013,Gratz2004}. To make matters worse, the increment of the international tourist trade and human mobility is promoting the spread of the main transmission vectors worldwide. For example, the {\it Aedes albopictus} has successfully invaded and occupied Africa, America, and even some cold areas of Europe \cite{Bhatt2013,Gratz2004}. Currently, there are no effective vaccines or therapeutic drugs to combat many mosquito-borne diseases. Thus, the primary strategy to prevent the spread of these diseases is to control the population of mosquitoes. Popular control strategies include the breeding places reduction, chemical insecticides, and the sterile insect technique (see \cite{abad2015mosquito,Dyck2021} and the references therein). This last strategy consists of releasing sterile mosquitoes (preferably male) so that a wild female mosquito that mates with a sterile male either does not reproduce or its produced eggs will not hatch. Sterile mosquitoes are not normally self-replicating and, therefore, cannot become established in the control region. The sterile insect technique has emerged recently as a biologically secure and effective tool for suppressing wild mosquitoes, with a notable performance in real-life problems, see \cite{Dyck2021,Lees2015} and the references therein.
	
	Understanding the interaction between wild and sterile mosquitoes is a critical step to improving the performance of any suppression strategy.
    A natural model for this biological situation is
    \begin{equation}\label{intro1}
		\left\{\begin{array}{lll}
			w'=w\left(\frac{a w}{w+g}-(\mu_{1}+\xi_{1}(w+g))\right),\\[4pt]
			g'=B-(\mu_{2}+\xi_{2}(w+g))g
		\end{array}\right.
	\end{equation}
    where $w$ and $g$ are the population densities of wild and sterile mosquitoes at time $t$, respectively, (see, for instance, \cite{Cai2014} and the references therein). In (\ref{intro1}), $a$ denotes the birth rate of wild mosquitoes, $\frac{w}{w+g}$ is the probability of a female mating with a non-sterile male, $B$ is the rate of releases of sterile mosquitoes, $\mu_{i}$ stands for the density-independent death rate, and 
	$\xi_{i}(w+g)$ represents the density-dependent self-regulation (density-dependent mortality due to intraspecific competition). In many situations, the sterile mosquitoes released do not have any influence on the population dynamics after losing their mating ability. Following this biological insight, the number of sterile mosquitoes released as an input nonnegative function instead of an independent variable, (see \cite{Yu2018} for more details on this modeling perspective). Employing this  idea in (\ref{intro1}), we arrive at 
	\begin{equation}\label{intro2}
		w'=w\left(\frac{a w}{w+g(t)}-(\mu_{1}+\xi_{1}(w+g(t)))\right)  
	\end{equation}
	where $g(t)$ is just a given function, typically piece-wise constant. Generally speaking, the model exhibits four dynamical scenarios, (see \cite{Wang2024, Yu2020a,Yu2020,Zheng2022,Zheng2023,Zheng2021b}):
	\begin{itemize}
		\item The trivial solution is a global attractor.
		\item There exists a non-trivial $T$-periodic solution which is a global attractor.
		\item There is a bi-stability between the origin and a non-trivial $T$-periodic solution.
		\item The trivial solution is a local attractor and a non-trivial $T$-periodic solution is semi-stable.
	\end{itemize}
	
	Despite its undoubted utility,  equation (\ref{intro2}) has an important limitation from a biological point of view:  the mortality and birth rates are not subject to environmental variations, an oversimplifying assumption in nature, (see \cite{bellver2024dynamics} and the references therein).  Many experimental works have indicated that the mosquitoes' growth is strongly influenced by humidity, daylight exposure, and temperature. Specifically, mosquitoes survive easily in humid ecosystems such as those with thick mat layers in the soil, (see \cite{Markos1958,Mulatti2014,Marini2016,Campos2012} and the references therein). In addition, mosquitoes' vital activity occurs mainly when there is little sunlight. For example, \cite{Mulatti2014,Marini2016} state that the main factors affecting the population size of \textit{Culex pipiens} (outside of intraspecific competition) are daylight length and temperature, both seasonally dependent. Higher temperatures in early spring can hasten the breeding season, while high temperatures during summer can be detrimental due to higher adult mortality.  Other studies such as \cite{Campos2012} show that humidity, precipitation and temperature contribute positively to the abundance of \textit{Aedes aegypti} in Brazil. On the other hand, it is broadly documented that the advent of harsh winters or dry seasons leads to a state of low metabolic activity, reduced morphogenesis, or limited physical activity in most mosquitoes, (see \cite{Hidalgo2018,Huestis2012,Denlinger2014} and the references therein). This phenomenon, known as diapause or aestivation depending on whether it is seasonal or it is triggered by external conditions, has a deep repercussion on the persistence of the population (see \cite{El-morshedy2024,Lou2019,Ruiz-herrera2022}).
	A possible extension of (\ref{intro2}) that solves the previous objections is
	\begin{equation}\label{intro3}
		w'=w\left(\frac{a(t) w}{w+g(t)}-(\mu_{1}(t)+\xi_{1}(t)(w+g(t)))\right).  
	\end{equation}
	It is worth mentioning that the introduction of time-dependent parameters makes the analysis considerably more difficult. In particular, the approach developed  in \cite{Wang2024, Yu2020a,Yu2020,Zheng2022,Zheng2023,Zheng2021b} is not valid for the nonautonomous counterpart (\ref{intro3}), and new ideas are needed.

In a series of papers \cite{duenas2023bifurcation,duenas2023critical,duenas2024critical,duenas2024generalized}, Dueñas, Núñez, and Obaya have developed a powerful approach to describe the dynamical behavior of scalar equations with d-concave nonlinearities, a property normally related to the sign of the third derivative. Motivated by these works, we show that for periodic equations coming from population dynamics, this type of condition is generally assocaited with  the positive sign of the third derivative of the inverse of the Poincaré map. This remark allows us to simplify some arguments in \cite{duenas2023bifurcation,duenas2023critical,duenas2024critical,duenas2024generalized} and obtain a deep geometrical understanding of the bifurcation patterns. Generally speaking, some degenerate situations cannot occur and the common transversality conditions at the bifurcation points are not required. It is worth mentioning that the results in \cite{duenas2023bifurcation,duenas2023critical,duenas2024critical,duenas2024generalized} are valid for a more general time-dependence, beyond periodicity. On the other hand, a contribution of this paper is to generalize the results in \cite{Wang2024, Yu2020a,Yu2020,Zheng2022,Zheng2023,Zheng2021b}
incorporating seasonality. Note that the arguments in those papers cannot be extended to time-dependent parameters.

	The structure of the paper is as follows. In Section 2, we prove some basic properties of model (\ref{intro3}) and study the sign of the third derivative of the inverse of its Poincaré map. From this analysis, we can deduce a first estimate of the number of periodic solutions of (\ref{intro3}). In Section 3, we describe the dynamical behavior of the model. The stability results and the bifurcation patterns are written in terms of the behavior of the origin, a crucial fact in model (\ref{intro3}). The reader can consult \cite{duenas2023bifurcation,duenas2023critical,duenas2024critical,duenas2024generalized,elia2025global} for similar global bifurcation patterns in nonautonomous scalar equations.  In Section 4, we apply our approach to other models considered in relevant literature (see \cite{Wang2024, Yu2020a,Yu2020,Zheng2022,Zheng2023,Zheng2021b} and the references therein). Finally, in the last section, we discuss our findings. One of the main conclusions of this paper is that seasonality is paramount for the performance of the sterile insect technique.

	\section{Estimating the number of $T$-periodic solutions}
	Consider the equation
	\begin{equation}\label{eq1} 
		w'=w\left(\frac{a(t)w}{w+g(t)}-\mu(t)-\xi(t)(w+g(t))\right)
	\end{equation}
	where the functions $a,\mu,\xi:\mathbb{R}\longrightarrow (0,+\infty)$ are of class $\mathcal{C}^{\infty}$ and $T$-periodic. We assume that $g:\mathbb{R}\longrightarrow [0,+\infty)$ is  $T$-periodic and piecewise constant of the form
	\begin{equation}\label{g}
			g(t)= \left\{\begin{array}{ll}
				g_{0} & \mathrm{if } t\in  [i T,iT + \overline{T}),
				\\[4pt]
				0& \mathrm{if } t\in [iT + \overline{T},(i+1)T),
			\end{array}\right.
			\quad 
			i\in \mathbb{Z},
	\end{equation}
	with $0<\overline{T}<T$ and $g_{0}>0$. Along this paper, $F:\mathbb{R}\times [0,+\infty)\longrightarrow \mathbb{R}$ denotes the map associated with (\ref{eq1}), that is,
	\begin{equation}\label{definitionF}
		F(t,w)=w\left(\frac{a(t)w}{w+g(t)}-\mu(t)-\xi(t)(w+g(t))\right).
	\end{equation}
	 All parameters of (\ref{eq1}) share the same period. This condition will be relaxed in Section 4.

	The discontinuity of the function $g$ allows us to visualize (\ref{eq1}) as a switching between two equations. Specifically, if $w(t;w_{0})$ denotes the solution of (\ref{eq1}) with initial condition $w_{0}\in [0,+\infty)$  defined on
	the maximal (right) interval $[0,\alpha)$, then $w(t;w_{0})$ is the solution of 
	\begin{equation}\label{eq2}
		\left\{\begin{array}{lll}
			w'=w\left(\frac{a(t)w}{w+g_{0}}-\mu(t)-\xi(t)(w+g_{0})\right)\\[4pt]
			w(0)=w_{0}
		\end{array}\right.
	\end{equation}
	in the interval $[0,\overline{T})\cap [0,\alpha)$, the solution of 
	\begin{equation}\label{eq3}
		\left\{\begin{array}{lll}
			w'=w\left(a(t)-\mu(t)-\xi(t)w\right)\\[4pt]
			w(\overline{T})=w(\overline{T};w_{0})
		\end{array}\right.
	\end{equation}
	in the interval $[\overline{T},T)\cap [0,\alpha)$; and the same procedure is $T$-periodically repeated in the other intervals. It is worth mentioning that $w(t;w_{0})$ is continuous in $[0,\alpha)$ and of class $\mathcal{C}^{\infty}$ in the intervals $[0,\alpha)\cap (iT,i T+\overline{T})$ and $[0,\alpha)\cap (i T+\overline{T},(i+1)T)$ for $i=0,1,2,\ldots$ Moreover, the one-sided derivatives of $w(t;w_{0})$ at $t\in \{i T+\overline{T}, (i+1) T:i=0,1,2,\ldots\} $ are finite.
	\begin{remark}
		The solution $w(t;w_{0})$ does not depend on the values of $g$ at the points  $t\in \{i T+\overline{T}, (i+1) T:i=0,1,2,\ldots\} $. For example, the function
		
		\begin{equation*}
			g(t)= \left\{\begin{array}{ll}
				g_{0} & \mathrm{if } t\in  (i T,iT + \overline{T}),
				\\[4pt]
				0& \mathrm{if } t\in [iT + \overline{T},(i+1)T],
			\end{array}\right.
			\quad 
			i\in \mathbb{Z},
		\end{equation*}
		
		 produces exactly the same solutions.
	\end{remark}
	Before proving that equation (\ref{eq1}) is well-posed, we state an auxiliary result. The proof is given in an Appendix for the reader's convenience.
	\begin{lemma}\label{l1}
		Let $G:\mathbb{R}\times[0,+\infty)\longrightarrow \mathbb{R}$ be a map of class $\mathcal{C}^{\infty}$ so that
		\begin{equation}\label{desi1}
			F(t,w)<G(t,w)  
		\end{equation}
		for all $(t,w)\in [0,+\infty)\times (0,+\infty)$. 
		Denote by $\widetilde{w}(t;w_{0})$ the solution of
		\begin{equation*}
			w'=G(t,w)
		\end{equation*}
		that satisfies $\widetilde{w}(0;w_{0})=w_{0}$. Given $w_{0}>0$, if $\widetilde{w}(t;w_{0})$ is defined in $[0,\beta)$ with $\beta>0$, then $w(t;w_{0})$ is defined in $[0,\beta)$ and
		\begin{equation*}
			w(t;w_{0})<\widetilde{w}(t;w_{0})
		\end{equation*}
		for all $t\in (0,\beta)$.
	\end{lemma}
	\begin{proposition}\label{p1}
		For any $w_{0}\in [0,+\infty)$, $w(t;w_{0})$ is defined  on $[0,+\infty)$, $w(t;w_{0})\geq 0$ for all $t\geq 0$, and  
		\begin{equation}\label{limite}
			\limsup_{t\rightarrow +\infty}w(t;w_{0})\leq \Gamma   
		\end{equation}
		with $\Gamma=\max\left\{\frac{a(t)}{\xi(t)}:t\in[0,T]\right\}.$
	\end{proposition}
	\begin{proof}
		Take an initial condition $w_{0}\in [0,+\infty)$. It is not restrictive to assume that $w_{0}>0$; otherwise $w(t;w_{0})=0$ and the conclusion is clear. We stress that $w(t;w_{0})>0$ for all $t\in [0,\alpha)$ by uniqueness of solution. Using that $\frac{a(t)w}{w+g_{0}}< a(t)$ for all $w\geq 0$, we have that 
		\begin{equation*}
			F(t,w)<G(t,w)
		\end{equation*}
		with $G(t,w)=w(a(t)-\xi(t) w)$ for all $(t,w)\in [0,+\infty)\times (0,+\infty)$. After the simple integration of the equation
		\begin{equation*}
			w'=G(t,w),
		\end{equation*}
		the solution 
		with initial condition $w_{0}$, say $\widetilde{w}(t;w_{0})$, is defined for all $t\geq 0$ and
		\begin{equation}
			\limsup_{t\to +\infty}\widetilde{w}(t;w_{0})\leq \Gamma.
		\end{equation}
		Now, the conclusion follows directly from Lemma \ref{l1}.
	\end{proof}
	
	\noindent As a direct consequence of the previous proposition,
	\begin{equation*}
		[0,+\infty)=\{w_{0}\in [0,+\infty):w(t;w_{0}) \mathrm{ is defined for all } t\in [0,T]\}.
	\end{equation*}
	In this framework, the dynamical behavior of (\ref{eq1}) is determined by the discrete equation
	\begin{equation*}
		x_{n+1}=P(x_{n})
	\end{equation*}
	with $P:[0,+\infty)\longrightarrow P([0,+\infty))$ the Poincaré map of (\ref{eq1}), that is,  $P(w_{0})=w(T;w_{0})$.
	For the reader's convenience, we highlight some useful properties of $P$, see \cite{Ortega2019} for more details.
	\begin{description}
		\item[(A1)] $P$ is strictly increasing and continuous. Therefore, $P([0,+\infty))$ is an open set (relative to $[0,+\infty)$).
		\item[(A2)] The fixed points of $P$ are the initial conditions that lead to $T$-periodic solutions of (\ref{eq1}).
		\item[(A3)] $P^{n}(w_{0})=P\underbrace{\circ\ldots\circ}_{n} P(w_{0})=w(nT;w_{0})$. Moreover, $P^{-1}:P([0,+\infty))\longrightarrow [0,+\infty)$ is given by $P^{-1}(w_{0})=w(-T;w_{0})$.
		\item[(A4)] $\lim_{n\to+\infty}P^{n}(w_{0})=q$ with $q\in [0,+\infty)$ is equivalent to $\lim_{t\to+\infty}[w(t;w_{0})-w(t;q)]=0$
		with $w(t;q)$ a $T$-periodic solution of (\ref{eq1}).
	\end{description}
	As mentioned above, (\ref{eq1}) can be visualized as a switching between equations (\ref{eq2}) and (\ref{eq3}). This ensures that $P$ is of class $\mathcal{C}^{\infty}$ in $[0,+\infty)$. Note that  $P=P_{2}\circ P_{1}$ with $P_{1}(w_{0})=w_{1}(\overline{T};w_{0})$ and $P_{2}(w_{0})=w_{2}(T;\overline{T},w_{0})$
	being $w_{1}(t;w_{0})$ the solution of 
	\begin{equation*}
		\left\{\begin{array}{lll}
			w'=w\left(\frac{a(t)w}{w+g_{0}}-\mu(t)-\xi(t)(w+g_{0})\right)\\[4pt]
			w(0)=w_{0}
		\end{array}\right.
	\end{equation*}
	and $w_{2}(t;\overline{T},w_{0})$ the solution of 
	\begin{equation*}
		\left\{\begin{array}{lll}
			w'=w\left(a(t)-\mu(t)-\xi(t)w\right)\\[4pt]
			w(\overline{T})=w_{0}.
		\end{array}\right.
	\end{equation*}
	Properties {\bf (A2)}-{\bf (A4)} and the continuity of $P$ were deduced in Section 1.2 of \cite{Ortega2019} for general systems in $\mathbb{R}^d$, {\it i.e.}, 
	\begin{equation*}
		x'=X(t,x)
	\end{equation*}
	with $X:\mathbb{R}\times \mathbb{R}^d\longrightarrow \mathbb{R}^d$  continuous, $T$-periodic, and the initial value problem having uniqueness of solution. We stress that exactly the same proofs in \cite{Ortega2019} are valid when we consider $T$-periodic switching systems like (\ref{eq1}). On the other hand, we observe that $P$ is strictly increasing as a direct consequence of the uniqueness of solution for the initial value problem  (\ref{eq1}).

	Let $\Fix(P)$ be the fixed point set of $P$. By Proposition \ref{p1}, $\Fix(P)$ is a bounded set. We employ the notation $\Delta=\max\{\Fix(P)\}$ in the sequel.

	\begin{proposition}\label{p2}
		$P(w_{0})<w_{0}$ for all $w_{0}\in (\Delta,+\infty)$.
	\end{proposition}
	\begin{proof}
		Using that $P$ is continuous and $P(w_{0})\not=w_{0}$ for all $w_{0}\in (\Delta,+\infty)$, we have that either $P(w_{0})>w_{0}$ for all $w_{0}\in (\Delta,+\infty)$ or 
		$P(w_{0})<w_{0}$ for all $w_{0}\in (\Delta,+\infty)$. Assume, by contradiction, that the first case holds. Then, by {\bf (A1)} and {\bf (A3)}, $P^{n}(w_{0})=w(nT;w_{0})$ is a strictly increasing sequence for all $w_{0}\in (\Delta,+\infty)$. Using that $\Fix(P)\cap (\Delta,+\infty)=\emptyset$, we conclude that $\{P^{n}(w_{0})\}=\{w(nT;w_{0})\}$ is unbounded. This is a contradiction with Proposition \ref{p2}.
	\end{proof}
	Next, we prove that $P$ has, at most, three fixed points. The main ingredients to deduce this result are the formulas for the successive derivatives of the Poincaré map obtained by Lloyd, (see pages 284-285 in \cite{Lloyd1979}).
	\begin{proposition}\label{p3}
		Consider the equation 
		\begin{equation}\label{aux1}
			x'=S(t,x)
		\end{equation}
		with $S:O\subset \mathbb{R}^{2}\longrightarrow \mathbb{R}$ a map of class $\mathcal{C}^{\infty}$ in $x$ and continuous in $t$. Assume that the solution $x(t;x_{0})$ of (\ref{aux1}) is defined in $[0,T]$. Then,
		\begin{equation*}
			P'(x_{0})=\exp\left(\int_{0}^{T}\frac{\partial S}{\partial x}(t,x(t;x_{0}))d{t}\right),
		\end{equation*}
		\begin{equation*}
			P''(x_{0})=P'(x_{0})\left(\int_{0}^{T}\frac{\partial^{2}S}{\partial x^{2}}(t,x(t;x_{0})) \exp\left(\int_{0}^{t}\frac{\partial S}{\partial x}(s,x(s;x_{0}))d{s}\right)d{t}\right),
		\end{equation*}
		and 
		\begin{equation*}
			P'''(x_{0})=P'(x_{0})\left(\frac{3}{2}\left(\frac{P''(x_{0})}{P'(x_{0})}\right)^2+\int_{0}^{T}\frac{\partial^3 S}{\partial x^3}(t,x(t;x_{0}))\exp\left(2\int_{0}^{t}\frac{\partial S}{\partial x}(s,x(s;x_{0}))d{s}\right)d{t}\right).
		\end{equation*}
	\end{proposition}
	\begin{theorem}\label{t1} Equation (\ref{eq1}) has, at most, three $T$-periodic solutions.
	\end{theorem}
	\begin{proof}
		First, we make the change of variable $s=-t$. For  $y(t)=w(-t)$, we obtain that
		\begin{equation}\label{change}
			y'=-y \left(\frac{\widetilde{a}(t) y}{y+\widetilde{g}(t)}-\widetilde{\mu}(t)-\widetilde{\xi}(t)(y+\widetilde{g}(t))\right)
		\end{equation}
		with $\widetilde{a}(t)=a(-t)$, $\widetilde{\mu}(t)=\mu(-t)$, $\widetilde{\xi}(t)=\xi(-t)$, and 
		
			\begin{equation*}
				\widetilde{g}(t)=g(-t)=\left\{\begin{array}{ll}
					0 & \text{if } t\in (-(i+1)T,-iT-\overline{T}],\\[4pt]
					g_{0} & \text{if } t\in (-iT-\overline{T},-iT]
				\end{array}\right.
			\end{equation*}
			with $i\in\mathbb{Z}$.
		It is clear that the number of $T$-periodic solutions of (\ref{eq1}) and (\ref{change}) is the same. Actually, if $\widetilde{P}$ is the Poincaré map of (\ref{change}), $\widetilde{P}=P^{-1}$, that is, the inverse of the Poincaré map of (\ref{eq1}). Thus, 
		\begin{equation*}
			\Fix(P)=\Fix(\widetilde{P})\subset [0,\Delta].
		\end{equation*}
		Since $P([0,+\infty))$ is an open set (relative to $[0,+\infty)$), we can take $K>0$ so that $\Delta<K$ and $[0,K]\subset P([0,+\infty))$. Our task now is to prove that $\widetilde{P}'''(w_{0})>0$ for all $w_{0}\in [0,K]$. We emphasize that this property ensures that $\widetilde{P}$ has, at most, three fixed points in $[0,K]$. Let
		\begin{equation*}
			H(t,y)=-y \left(\frac{\widetilde{a}(t)y}{y+\widetilde{g}(t)}-\widetilde{\mu}(t)-\widetilde{\xi}(t)(y+\widetilde{g}(t))\right).
		\end{equation*}
		Using that $\frac{\partial^{3} H}{\partial y^3}(t,y)=\frac{6 \widetilde{a}(t)\widetilde{g}(t)^2}{(y+\widetilde{g}(t))^4}$ and Proposition \ref{p3}, we conclude that the third derivative of $\widetilde{P}$ is strictly positive in $[0,K]$.

A possible objection of the previous argument might be that $H$ is not continuous at $t$. Next, we prove that the formulas in Proposition \ref{p3} are also valid for (\ref{change}). Indeed,
	$$y(t;y_{0})=\left\{\begin{array}{lll}
		y_{0}+\int_{0}^{t}H(s;y(s;y_{0}))ds\;\;\;\;\;\;\;\;\;\;\;\;\;\;\;\;\;\;\;\;\;\;\;\;\;\;\;\;\;\;\;\;\;\;\;\;\;\;\;\;\;\;\;\;\;\;\;\;\;{\rm if}\;\;t\in [0,T-\overline{T}]\\[4pt]
		y_{0}+\int_{0}^{T-\overline{T}}H(s;y(s;y_{0}))ds+\int_{T-\overline{T}}^{t}H(s;y(s;y_{0}))ds\;\;\;\;\;\;{\rm if}\;\;t\in [T-\overline{T},T].
	\end{array}\right.$$
	Since the integral terms $\int_{0}^{t}H(s;y(s;y_{0}))ds$, $\int_{0}^{T-\overline{T}}H(s;y(s;y_{0}))ds$, and $\int_{T-\overline{T}}^{t}H(s;y(s;y_{0}))ds$ are of class $\mathcal{C}^{\infty}$ (with respect to $y_{0}$), we have that
	$$\frac{\partial y(t;y_{0})}{\partial y_{0}}=\left\{\begin{array}{lll}
		1+\int_{0}^{t}\frac{\partial}{\partial y}H(s;y(s;y_{0}))\frac{\partial y(s;y_{0})}{\partial y_{0}}ds\;\;\;\;\;\;\;\;\;\;\;\;\;\;\;\;\;\;\;\;\;\;\;\;\;\;\;\;\;\;\;\;\;\;\;\;\;\;\;\;\;\;\;\;\;\;\;\;\;{\rm if}\;\;t\in [0,T-\overline{T}]\\[4pt]
		1+\int_{0}^{T-\overline{T}}\frac{\partial}{\partial y}H(s;y(s;y_{0})) \frac{\partial y(s;y_{0})}{\partial y_{0}}ds+\int_{T-\overline{T}}^{t}\frac{\partial}{\partial y}H(s;y(s;y_{0})) \frac{\partial y(s;y_{0})}{\partial y_{0}}ds\;\;\;\;\;\;{\rm if}\;\;t\in [T-\overline{T},T].
	\end{array}\right.$$
	Then, $\frac{\partial y(t;y_{0})}{\partial y_{0}}$ is solution of the problem
	$$\left\{\begin{array}{lll}z'=\frac{\partial H}{\partial y}(t;y(t;y_{0}))z\\z(0)=1.\end{array}\right.$$
	
	Consequently, $\widetilde{P}'(y_{0})=\frac{\partial y(T;y_{0})}{\partial y_{0}}=e^{\int_{0}^{T}\frac{\partial H(s;y(s;y_{0}))}{\partial y}ds}$. This argument proves that the formula for the first derivative in Proposition \ref{p3} is valid for (\ref{change}). To deduce the formulas for $\widetilde{P}''$ and $\widetilde{P}'''$, we have to repeat the same argument as that in \cite{Lloyd1979} splitting the integral terms as above.

	\end{proof}
	\begin{remark}
	 Following Dueñas {\it et al.}  \cite{duenas2023bifurcation,duenas2023critical,duenas2024critical,duenas2024generalized}, model (\ref{eq1}) is d-concave. A similar theorem was derived in \cite{duenas2023bifurcation} with a different approach, see Section 4. For the periodic case, the use of Proposition 2.3 simplifies considerably the arguments.
	\end{remark}

	\section{Global dynamics of model (\ref{eq1})}
	In the previous section, we have proved that (\ref{eq1}) has, at most, three $T$-periodic solutions, {\it i.e.}, the trivial solution (that always exists) and, at most, two non-trivial $T$-periodic solutions. Next, we analyze the dynamical behavior of model (\ref{eq1}). We approach this task through a deep analysis of the Poincaré map, a different methodology from that in  \cite{duenas2023bifurcation,duenas2023critical,duenas2024critical,duenas2024generalized}. On the other hand, despite a considerable amount of literature on this model (see \cite{Wang2024, Yu2020a,Yu2020,Zheng2022,Zheng2023,Zheng2021b} and the references therein), to the best of our knowledge, the analysis of model (\ref{eq1}) with time-dependent parameters has not been provided in the literature yet.

    To avoid misleading conclusions, we state some basic notions of stability theory. We say that a $T$-periodic solution $w(t;w_{0}^*)$ of (\ref{eq1}) is asymptotically stable (resp. unstable) if there is an open neighborhood $U$ (relative to $[0,+\infty)$) of $w_{0}^*$ so that $\lim_{n\rightarrow +\infty}P^{n}(w_{0})=w_{0}^*$ (resp. $\lim_{n\rightarrow +\infty}P^{-n}(w_{0})=w_{0}^*$) for all $w_{0}\in U$ with $P$ the Poincaré map of (\ref{eq1}). We recall that $P'(w_{0}^*)<1$ (resp. $P'(w_{0}^*)>1$) implies that $w(t;w_{0}^*)$ is asymptotically stable (resp.  unstable). We say that a $T$-periodic solution $w(t;w_{0}^*)$ is an attractor in $V\subset[0,+\infty)$ if $\lim_{n\rightarrow +\infty} P^{n}(x)=w_{0}^*$ for all $x\in V$. 
	As emphasized in property {\bf (A4)}, these notions can be formulated regarding the solutions. For example, $w(t;w_{0}^*)$  is asymptotically stable  if there is an open neighborhood $U$ (relative to $[0,+\infty)$) of $w_{0}^*$ so that $\lim_{t\rightarrow +\infty}[w(t;w_{0})-w(t;w_{0}^*)]=0$ for all $w_{0}\in U$. 
	
	The next result characterizes the stability of the trivial solution. For simplicity in the exposition, we employ the notation
	\begin{equation}\label{i1}
		I_{1}=-\int_{0}^{\overline{T}}\xi(t)g_{0}d{t}-\int_{0}^{T}\mu(t)d{t}+\int_{\overline{T}}^{T}a(t)d{t} .
	\end{equation}
	
	\begin{proposition}\label{p4}
		\begin{itemize}
			\item[i)] If $I_1 > 0$, then $P'(0)>1$ and  the origin is unstable.
			
			\item[ii)] If $I_1 < 0$, then $P'(0)<1$ and the origin is asymptotically stable.
		\end{itemize}
	\end{proposition}
	\begin{proof}
		First, we derive a useful representation of the Poincaré map. Given $w_{0}\in (0,+\infty)$, the solution $w(t;w_{0})$ satisfies 
		\begin{equation*}
			w'(t;w_{0})=w(t;w_{0})\left(\frac{a(t)w(t;w_{0})}{w(t;w_{0})+g(t)}-\mu(t)-\xi(t)(w(t;w_{0})+g(t))\right),
		\end{equation*}
		or equivalently, 
		\begin{equation*}
			\frac{w'(t;w_{0})}{w(t;w_{0})}=\frac{a(t)w(t;w_{0})}{w(t;w_{0})+g(t)}-\mu(t)-\xi(t)(w(t;w_{0})+g(t)).
		\end{equation*}
		Integrating over $0$ and $T$, we deduce that
		\begin{equation*}
			\ln w(T;w_{0})-\ln w_{0}=\int_{0}^{T} \paren{\frac{a(s) w(s;w_{0})}{w(s;w_{0})+g(s)}-\mu(s)-\xi(s)(w(s;w_{0})+g(s))} d{s}.
		\end{equation*}
		This implies that
		\begin{align}\label{expresionpoin}
			& P(w_{0}) = w(T;w_{0})=w_{0}\exp\left(\int_{0}^{T} \paren{\frac{a(s) w(s;w_{0})}{w(s;w_{0})+g(s)}-\mu(s)-\xi(s)(w(s;w_{0})+g(s))} d{s}\right)
			\\ \notag & = w_{0}\exp\left(\int_{0}^{\overline{T}} \paren{\frac{a(s) w(s;w_{0})}{w(s;w_{0})+g_{0}}-\mu(s)-\xi(s)(w(s;w_{0})+g_{0})} d{s}+\int_{\overline{T}}^{T} \paren{a(s)-\mu(s)-\xi(s)w(s;w_{0})} d{s}\right).
		\end{align}
		Thus, $P(w_{0})=w_{0}\Psi(w_{0})$ with 
		\begin{equation*}
			\Psi(w_{0})=\exp\left(\int_{0}^{\overline{T}} \paren{\frac{a(s) w(s;w_{0})}{w(s;w_{0})+g_{0}}-\mu(s)-\xi(s)(w(s;w_{0})+g_{0})} d{s}+\int_{\overline{T}}^{T} \paren{a(s)-\mu(s)-\xi(s)w(s;w_{0})} d{s}\right).
		\end{equation*}
		Note that $P'(0)=\Psi(0)$ by definition since 
		\begin{equation*}
			P'(0) = \lim_{x \to 0^+} \frac{P(x) - P(0)}{x} = \lim_{x \to 0^+} \frac{P(x)}{x} = \lim_{x \to 0^+} \Psi(x) = \Psi(0) .
		\end{equation*}
		Using that $w(s;0)=0$, we directly conclude the thesis of the proposition.
	\end{proof}
	The following result states that the instability of the trivial solution guarantees the existence of a non-trivial $T$-periodic solution that is a global attractor in $(0,+\infty)$.
	
\begin{figure}[H]
		\centering
		\subfigure[Situation 1]{
			\includegraphics[width=0.475\linewidth]{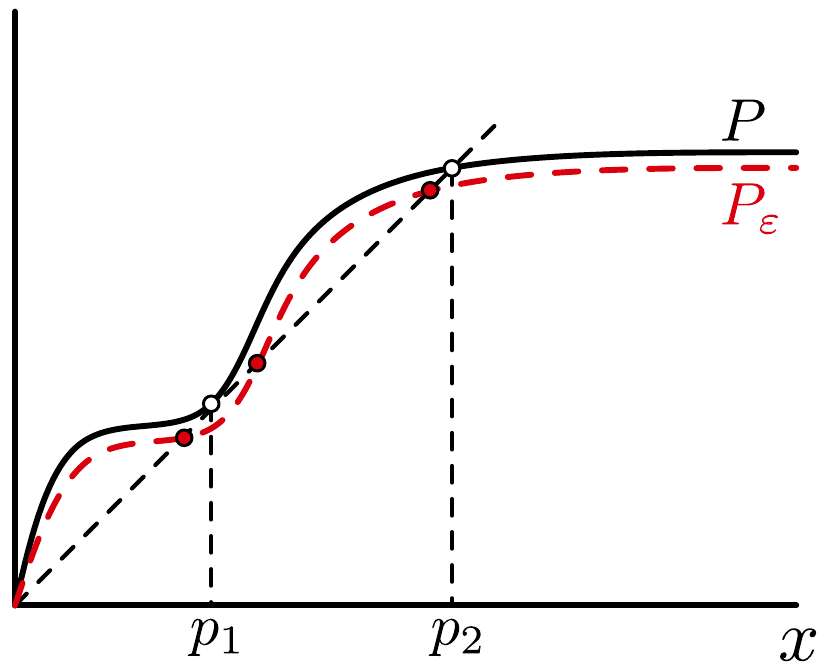}
		}
		\subfigure[Situation 2]{
			\includegraphics[width=0.475\linewidth]{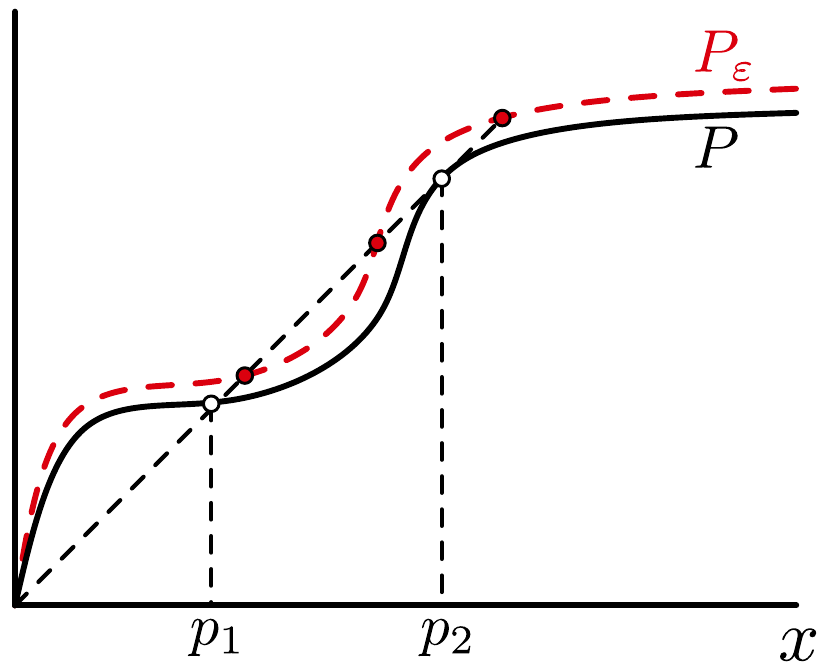}
		}
		\caption{Pictorial illustration of the proof of Theorem \ref{t2}}
		\label{fig:CasosTeorema3.1}
	\end{figure}

	\begin{theorem}\label{t2} 
		If \begin{equation}\label{inestabilidad}
			I_1 = -\int_{0}^{\overline{T}}\xi(t)g_{0}d{t}-\int_{0}^{T}\mu(t)d{t}+\int_{\overline{T}}^{T}a(t)d{t}>0,
		\end{equation}
		then there exists a non-trivial $T$-periodic solution $w(t;w_{0}^*)$ of model (\ref{eq1}) that is an attractor in $(0,+\infty)$, i.e., for all $w_{0}\in (0,+\infty)$,
		\begin{equation*}
			\lim_{t\rightarrow +\infty}[w(t;w_{0})-w(t;w_{0}^*)]=0.
		\end{equation*}
	\end{theorem}
	\begin{proof}
		We divide the proof into three steps:\newline
		\noindent {\bf Step 1:} Existence of a non-trivial $T$-periodic solution for equation (\ref{eq1}). \newline
		\noindent We know by Proposition \ref{p4} that $P'(0)>1$ by (\ref{inestabilidad}). Therefore, there is $\delta>0$ so that $P(x)>x$ for all $x\in (0,\delta)$. On the other hand, $P(x)<x$ for all $x\in (\Delta,+\infty)$
		with $\Delta=\max\{\Fix(P)\}$, (see Proposition \ref{p2}). Now, the existence of a positive fixed point of $P$ is clear.\newline\newline
		\noindent {\bf Step 2:} Equation (\ref{eq1}) has exactly one non-trivial $T$-periodic solution.\newline
		Assume, by contradiction, that (\ref{eq1}) has exactly two non-trivial $T$-periodic solutions. Let $p_{1},p_{2}$ be the positive fixed points of $P$ with $0<p_{1}<p_{2}$. Since $P'(0)>1$ (see Proposition \ref{p4}), we deduce that $P(x)>x$ for all $x\in (0,p_{1})$. On the other hand, $P(x)<x$
		for all $x\in (p_{2},+\infty)$ (by Proposition \ref{p2}). Then, there are two possible situations:
		\begin{itemize}
			\item {\bf Situation 1:} $P(x)>x$ if $x\in (p_{1},p_{2})$.
			\item {\bf Situation 2:} $P(x)<x$ if $x\in (p_{1},p_{2})$.
		\end{itemize}
		Now, we assume that situation 1 holds, (the differences regarding situation 2 are discussed at the end).
		First, we pick three points $q_0,q_{1},q_{2}$ with $q_{0}\in(0,p_{1})$, $q_{1}\in (p_{1},p_{2})$, and $q_{2}\in (p_{2},+\infty)$ . We have that 
		\begin{equation}\label{condt1}
			P(q_{0})=w(T,q_{0})>q_{0}, \mathrm{ } P(q_{1})= w(T;q_{1})>q_{1}, \mathrm{ and } P(q_{2})=w(T;q_{2})<q_{2} .
		\end{equation}
		Next, we take $\varepsilon>0$ small enough so that
		
		\begin{equation}\label{perturb2}
			w_{\varepsilon}(T,q_{0})>q_{0}, \mathrm{ } w_{\varepsilon}(T;q_{1})>q_{1}, \mathrm{ and } w_{\varepsilon}(T;q_{2})<q_{2}
		\end{equation}
		where $w_{\varepsilon}(t;w_{0})$ is the solution of 
		\begin{equation}\label{perturb3}
			w'=w\left(\frac{a(t) w}{w+g(t)}-(\mu(t)+\varepsilon)-\xi(t)(w+g(t))\right)
		\end{equation}
		with $w_{\varepsilon}(0;w_{0})=w_{0}$. We can find such an $\varepsilon>0$ by continuous dependence of (\ref{perturb3}) with respect to $\varepsilon$ and (\ref{condt1}). Let $P_{\varepsilon}:[0,+\infty)\longrightarrow P_{\varepsilon}([0,+\infty))$ be the Poincaré map associated with equation (\ref{perturb3}).
		By (\ref{perturb2}), 
		\begin{equation}\label{finalt2}
			P_{\varepsilon}(q_{0})>q_{0}, \mathrm{ } P_{\varepsilon}(q_{1})>q_{1}, \mathrm{ and } P_{\varepsilon}(q_{2})<q_{2} .
		\end{equation}
		Moreover, by Theorem \ref{t1}, we know that $P_{\varepsilon}$ has, at most, three fixed points. On the other hand, we note that 
		\begin{equation*}
			w_{\varepsilon}(t;p_{1})<w(t;p_{1})
		\end{equation*}
		for all $t\in(0,T]$ because $F_{\varepsilon}(t,w)<F(t,w)$ for all $(t,w)\in [0,T]\times (0,+\infty)$
		with $F_{\varepsilon}$ and $F$ the maps associated with (\ref{perturb3}) and (\ref{eq1}), respectively. To deduce this property, we have to employ the same argument as that in Lemma \ref{l1}. Thus, \begin{equation}\label{finalt3}
			w_{\varepsilon}(T;p_{1})=P_{\varepsilon}(p_{1})<w(T;p_{1})=P(p_{1})=p_{1}.
		\end{equation}
		Collecting  (\ref{finalt2}) and (\ref{finalt3}), we conclude that $P_{\varepsilon}$ has, at least, four fixed points: the origin; a fixed point in $(q_{0},p_{1})$; a fixed point in $(p_{1},q_{1})$; and a fixed point in $(q_{1},q_{2})$. This is a contradiction. See Figure \ref{fig:CasosTeorema3.1}(a) for a pictorial explanation of this argument.  A similar argument works when situation 2 holds. For the reader's convenience, we sketch the main differences. Pick $q_{0}<q_{1}<q_{2}$ with $q_{0}\in (0,p_{1})$, $q_{1}\in (p_{1},p_{2})$ and $q_{2}\in (p_{2},+\infty)$;
		and  $\varepsilon>0$ small enough so that
		
		\begin{equation}\label{perturb2'}
			w_{\varepsilon}(T;q_{0})>q_{0}, \mathrm{ }w_{\varepsilon}(T;q_{1})<q_{1}, \mathrm{ and } w_{\varepsilon}(T;q_{2})<q_{2},
		\end{equation}
		where $w_{\varepsilon}(t;w_{0})$ is the solution of 
		\begin{equation}\label{perturb3'}
			w'=w\left(\frac{a(t) w}{w+g(t)}-(\mu(t)-\varepsilon)-\xi(t)(w+g(t))\right)
		\end{equation}
		with $w_{\varepsilon}(0;w_{0})=w_{0}$. Next, we can prove that the Poincaré map $P_{\varepsilon}$ has, at least, four fixed points. See Figure \ref{fig:CasosTeorema3.1}(b) for an illustration of this argument. This would be a contradiction.\newline\newline
		\noindent {\bf Step 3: Conclusion.}\newline
		The Poincaré map has a unique positive fixed point, say $w_{0}^*$. Moreover, the origin is unstable (see (\ref{inestabilidad}) and Proposition \ref{p4}). We also have that the orbits of 
		\begin{equation*}
			x_{n+1}=P(x_{n})
		\end{equation*}
		are bounded (see Proposition \ref{p1}). Using that $P$ is strictly increasing, we conclude that any orbit with a strictly positive initial condition tends to a fixed point. Thus, for any $w_{0}\in (0,+\infty)$, $\lim_{n\to+\infty}P^{n}(w_{0})=w_{0}^*$ or, equivalently,
		$\lim_{t\rightarrow+\infty}[w(t;w_{0})-w(t;w_{0}^*)]=0.$
	\end{proof}
	
The same conclusion as in the previous theorem when $a,\mu,$ and $\xi$ are constants was derived in \cite{Yu2020a}. In other words, we extend the results in \cite{Yu2020a}   when seasonality is introduced into the vital parameters of the population.

	To finish this section, we offer a complete description of the dynamical behavior of model (\ref{eq1}) when the reverse inequality in (\ref{inestabilidad}) is satisfied.
	\begin{theorem}\label{t3}
		Assume
		\begin{equation}\label{condaux4}
			I_1 = -\int_{0}^{\overline{T}}\xi(t)g_{0}d{t}-\int_{0}^{T}\mu(t)d{t}+\int_{\overline{T}}^{T}a(t)d{t}<0.
		\end{equation}
		Then, one of the following cases is satisfied:
		\begin{description}
			\item[i)] Equation (\ref{eq1}) does not have non-trivial $T$-periodic solutions. Moreover, for each $w_{0}\in [0,+\infty)$, $\lim_{t\to+\infty}w(t;w_{0})=0$.
			\item[ii)] Equation (\ref{eq1}) has a unique non-trivial $T$-periodic solution, say $w(t;w_{0}^*)$. Moreover, for each $w_{0}<w_{0}^*$, $\lim_{t\to+\infty}w(t;w_{0})=0$ and for each $w_{0}>w_{0}^*$, $\lim_{t\to+\infty}[w(t;w_{0})-w(t;w_{0}^*)]=0$.
			\item[iii)] Equation (\ref{eq1}) has exactly two non-trivial $T$-periodic solutions, say $w(t;w_{1}^*)$ and $w(t;w_{2}^*)$ with $w_{1}^*<w_{2}^*$. Moreover, for each $w_{0}\in (0,w_{1}^*)$, $\lim_{t\to+\infty}w(t;w_{0})=0$ and for each $w_{0}\in (w_{1}^*,+\infty)$, $\lim_{t\to+\infty}[w(t;w_{0})-w(t;w_{2}^*)]=0.$
		\end{description}
	\end{theorem}
	\begin{proof}
		Using (\ref{condaux4}) and Proposition \ref{p2}, $P'(0)<1$. We also know by Theorem \ref{t1} that $P$ has, at most, three fixed points, the origin, and, at most, two positive fixed points. Next, we reason depending on the number of positive fixed points of $P$.\newline\newline
		\noindent {\bf Case 1:} $P$ has no positive fixed points.\newline
		In this case, $P(x)<x$ for all $x\in (0,+\infty)$. Using that  $P$ is strictly increasing and $P(0)=0$, we conclude that $\lim_{n\to +\infty}P^n(w_{0})=0$ for all $w_{0}\in (0,+\infty)$. By {\bf (A4)}, we arrive at {\bf i)}. \newline\newline
		\noindent {\bf Case 2:} $P$ has a unique positive fixed point.\newline
		Let $w_{0}^*$ be such a fixed point. By Proposition \ref{p2}, we know that $P(x)<x$ for all $x\in (w_{0}^*,+\infty)$. Moreover, $P'(0)<1$ and $P(0)=0$ lead to $P(x)<x$ for all $x\in(0,w_{0}^*)$. Using that $P$ is increasing, we obtain that $\lim_{n\to+\infty} P^{n}(w_{0})=w_{0}^*$ for all $w_{0}\in (w_{0}^*,+\infty)$ and $\lim_{n\to+\infty} P^{n}(w_{0})=0$ for all $w_{0}\in (0,w_{0}^*)$. By {\bf (A4)}, we deduce {\bf ii)}.\newline\newline
		\noindent{\bf Case 3:} $P$ has exactly two positive fixed points.\newline
		Let $w_{1}^*,w_{2}^*$ be the positive fixed points of $P$ with $0<w_{1}^*<w_{2}^*$. At this moment, using that $P'(0)<1$ and Proposition \ref{p2}, we have that $P(x)<x$ for all $x\in (0,w_{1}^*)\cup (w_{2}^*,+\infty)$. There are two possible sub-cases for $P$ in the interval $(w_{1}^*,w_{2}^*)$:
		\begin{description}
			\item[Sub-case 1] $P(x)>x$ if $x\in (w_{1}^*,w_{2}^*)$.
			\item[Sub-case 2] $P(x)<x$ if $x\in (w_{1}^*,w_{2}^*)$.
		\end{description}
		To discard the possibility of {\bf Sub-case 2}, we  employ the same argument as that of excluding  Situation 1 in the previous theorem. In this case, we  consider the equation
			\begin{equation*}
				w'=w\left(\frac{a(t)w}{w+g(t)}-(\mu(t)-\varepsilon)-\xi(t)(w+g(t))\right)
			\end{equation*}
			with $\varepsilon>0$ small enough. Now, we would have that $P_{\varepsilon}$ has, at least, five fixed points. This would be a contradiction.

		Collecting the above information, we conclude that 
		\begin{equation*}
			x_{n+1}=P(x_{n})
		\end{equation*}
		has three fixed points, namely, $0,w_{1}^*$ and $w_{2}^*$ with $0<w_{1}^*<w_{2}^*$. Moreover, we know that $P$ is increasing and $P(x)<x$ for all $x\in (0,w_{1}^*)\cup (w_{2}^*,+\infty)$ and $P(x)>x$ for all $x\in (w_{1}^*, w_{2}^*)$. These properties directly imply {\bf iii)} by {\bf (A4)}.
	\end{proof}
	
The previous result describes the dynamical picture of (\ref{eq1}) when the origin is asymptotically stable. The same conclusion was obtained in \cite{Yu2020} when $a,\mu,\xi$ are constants. We stress that his approach is not valid in our framework.

	As a direct application of the proof of the previous theorem, we can deduce the following.
	\begin{corollary}\label{c1}
		Assume (\ref{condaux4}). The following conditions are equivalent:
		\begin{description}
			\item[i)] Equation (\ref{eq1}) has exactly two non-trivial $T$-periodic solutions.
			\item[ii)] There is $K>0$ so that $P(K)>K$ with $P$ the Poincaré map of (\ref{eq1}).
		\end{description}
	\end{corollary}
	There are many manners to give sufficient conditions that guarantee condition {\bf ii)} of Corollary \ref{c1}. For example, if there is a positive constant $K>0$ so that
	\begin{equation}\label{e1}
		\frac{a(t)K}{K+g_{0}}-\mu(t)-\xi(t)(K+g_{0})>0
	\end{equation}
	for all $t\in[0,\overline{T}]$ and
	\begin{equation}\label{e2}
		a(t)-\mu(t)-\xi(t)K>0
	\end{equation}
	for all $t\in[\overline{T},T]$, condition {\bf ii)} in Corollary \ref{c1} is satisfied. We stress  that (\ref{e1}) and (\ref{e2}) say that $F(t,K)>0$ for all $t\in [0,T]$, (see (\ref{definitionF})). For the choice $a=2$, $\mu=1$, $\xi=0.5$, $T=1$, $\overline{T}=0.75$, $K=1$ and $g_{0}=0.1$, equation (\ref{eq1}) satisfies (\ref{condaux4}) and condition {\bf ii)} of Corollary \ref{c1}.

		The expected bifurcations at the origin occur when $I_{1}=0$. The next result classifies these bifurcations.
		\begin{proposition}\label{pbif}
			Assume that $I_{1}=0$.
			\begin{description}
				\item[i)] If $P''(0)>0$, then there exists a non-trivial $T$-periodic solution $w(t;w_{0}^*)$ of model (\ref{eq1}) that is an attractor in $(0,+\infty)$.
				\item[ii)] If $P''(0)\leq 0$, then $\lim_{t\to+\infty}w(t;w_{0})=0$ for each $w_{0}\in [0,+\infty)$.
			\end{description}
		\end{proposition}
		\begin{proof}
			{\bf i)} It is clear that $P'(0)=1$ and $P''(0)>0$ implies that $0$ is unstable. The conclusion is now a direct consequence of the proof of Theorem \ref{t2}.  Note that the unique role of (\ref{inestabilidad}) is to guarantee that the origin is unstable.\newline
			\noindent{\bf ii)} First, we observe that $P''(0)\leq 0$ implies that the origin is asymptotically stable. This is clear when $P''(0) < 0$. If $P''(0) = 0$,  we deduce $P'''(0) < 0$  by using Lloyd's formula. Then, the origin is also asymptotically stable because $P'(0)=1$, $P''(0)=0$, and $P'''(0)<0$. Next, we discard {\bf ii)} and {\bf iii)} in Theorem \ref{t3}. We stress that (\ref{condaux4}) is employed in Theorem \ref{t3} to deduce that the origin is asymptotically stable. Assume, by contradiction, that {\bf ii)} holds. Then, there is a positive fixed point $w_{0}^*>0$ so that $P(w)<w$ for all $w>0$ with $w\not=w_{0}^*$. Take two points $q_{1},q_{2}$ with $q_{1}\in (0,w_{0}^*)$ and $q_{2}\in (w_{0}^*,+\infty).$ We know that
			\begin{equation*}
				P(q_{1})<q_{1} \mathrm{ and }P(q_{2})<q_{2}.
			\end{equation*}
			Next, we pick $\varepsilon>0$ small enough so that
			\begin{equation*}
				P_{\varepsilon}(q_{1})<q_{1} \mathrm{ and }P_{\varepsilon}(q_{2})<q_{2}
			\end{equation*}
			with $P_{\varepsilon}$ the Poincaré map of
			\begin{equation*}
				w'=w\left(\frac{a(t) w}{w+g(t)}-(\mu(t)-\varepsilon)-\xi(t)(w+g(t))\right).
			\end{equation*}
			Using Proposition \ref{p4} and $I_{1}=0$, we deduce that $P_{\varepsilon}'(0)>1$. Collecting the previous information, we conclude that $P_{\varepsilon}$ has, at least, three positive fixed points. This is a contradiction with Theorem \ref{t1}. The argument for excluding {\bf iii)} is analogous.\end{proof}
		\color{black}
		After simple computations (see Appendix), we have that 
			\begin{align}
				\label{i1tilde}
				P''(0) &=  2\int_{0}^{\overline{T}} \paren{\frac{a(s)}{g_0} - \xi(s)} \exp\paren{-\int_{0}^{s} (g_0 \xi(t) + \mu(t)) d{t} } d{s} -
				\\[10pt]
				\label{i1tilde2}
				& - 2\int_{\overline{T}}^{T} \xi(s) \exp\paren{\int_{\overline{T}}^{s}a(t) d{t} - \int_{0}^{\overline{T}} g_0 \xi(t) d{t} - \int_{0}^{s} \mu(t) d{t} } d{s}.
		\end{align}
		\color{black}
		Theorem \ref{t3} classifies the dynamical behavior of (\ref{eq1}) when (\ref{condaux4}) holds. A common fact is that the origin is always a local attractor. Our next task consists of quantitatively estimating its basin of attraction. For simplicity, define
		\begin{equation}\label{i2}
			I_{2}=-\int_{0}^{\overline{T}}\xi(t)g_{0}d{t}-\int_{0}^{T}\mu(t)d{t}+\int_{0}^{T}a(t)d{t}.
		\end{equation}
		Using the definition of $I_{1}$ (see \eqref{i1}), we realize that
		\begin{equation*}
			I_{2}-I_{1}=\int_{0}^{\overline{T}}a(t)dt.
		\end{equation*}

		\begin{lemma}\label{lb1} If $I_{2}\leq 0$, $\lim_{t\rightarrow +\infty}w(t;w_{0})=0$ for all $w_{0}\in (0,+\infty)$.
		\end{lemma}
		\begin{proof}
			Consider the equation
			\begin{equation}\label{b1}
				w'=w(a(t)-\mu(t)-\xi(t)g(t)-\xi(t)w).
			\end{equation}
			After a simple integration of the equation, we know that $0$ is a global attractor of (\ref{b1}) provided $I_{2}\leq 0$. The conclusion now follows from Lemma \ref{l1}.
		\end{proof}
		
		\begin{proposition}\label{pb1} Assume that (\ref{condaux4}) holds together with the following conditions:
			\begin{description}
				\item[i)] $a(t)>\mu(t)$ for all $t\in [0,\overline{T}]$.
				\item[ii)] $I_{2}>0$.
			\end{description}
			Then, $\lim_{t\rightarrow +\infty}w(t;w_{0})=0$ for all $w_{0}\in \left(0,\frac{-I_{1}g_{0}}{I_{2}e^{\int_{0}^{\overline{T}} \paren{a(t)-\mu(t)} d{t}}}\right)$.
		\end{proposition}
		\begin{proof}
			We divide the proof into two steps.\newline
			\noindent {\bf Step 1:} $P(w_{0})<w_{0}$ for all $w_{0}\in  \left(0,\frac{-I_{1}g_{0}}{I_{2}e^{\int_{0}^{\overline{T}} \paren{a(t)-\mu(t)} d{t}}}\right) $.\newline\newline
			Using the expression of the Poincaré map derived in Proposition \ref{p4}, see (\ref{expresionpoin}), we have that $P(w_{0})<w_{0}$ if, and only if,
			\begin{equation*}
				\int_{0}^{\overline{T}} \paren{\frac{a(s) w(s;w_{0})}{w(s;w_{0})+g_{0}}-\mu(s)-\xi(s)(w(s;w_{0})+g_{0})} d{s}+\int_{\overline{T}}^{T} \paren{a(s)-\mu(s)-\xi(s)w(s;w_{0})} d{s}<0.
			\end{equation*}
			After a reordering of terms, $P(w_{0})<w_{0}$ is equivalent to
			\begin{equation}\label{b3}
				\int_{0}^{\overline{T}}\frac{a(t) w(t;w_{0})}{w(t;w_{0})+g_{0}}d{t}<-I_{1}+\int_{0}^{T}\xi(t)w(t;w_{0})d{t} .
			\end{equation}
			On the other hand, 
			\begin{equation*}
				w'(t;w_{0})\leq w(t;w_{0})(a(t)-\mu(t))
			\end{equation*}
			for all $t\geq 0$. This implies that $w(t;w_{0})\leq w_{0}e^{\int_{0}^{t}(a(s)-\mu(s))d{s}}$ for all $t\geq 0$. We note that $w_{0}e^{\int_{0}^{t} (a(s)-\mu(s))d{s}}\leq w_{0}e^{\int_{0}^{\overline{T}} (a(s)-\mu(s))d{s}}$ for all $t\in [0,\overline{T}]$ by {\bf i)}. Using that $h(x)=\frac{x}{x+g_{0}}$ is strictly increasing in $[0,+\infty)$, we conclude that (\ref{b3}) is satisfied when
			\begin{equation}\label{b6}
				\frac{w_{0}e^{\int_{0}^{\overline{T}}\paren{a(s)-\mu(s)} d{s}}}{g_{0}+w_{0}e^{\int_{0}^{\overline{T}}\paren{a(s)-\mu(s)} d{s}}}\int_{0}^{\overline{T}} a(t)d{t}<-I_{1}
			\end{equation}
			holds. By a straightforward computation,  
			\begin{equation*}
				w_{0}\in \left(0,\frac{-I_{1}g_{0}}{I_{2}e^{\int_{0}^{\overline{T}}\paren{a(t)-\mu(t)}d{t}}}\right)
			\end{equation*}
			implies (\ref{b6}).
			
			\noindent {\bf Step 2: Conclusion.}\newline
			\noindent We know that $P(0)=0$ and $P(w_{0})<w_{0}$ for all $w_{0}\in \left(0,\frac{-I_{1}g_{0}}{I_{2}e^{\int_{0}^{\overline{T}} \paren{a(t)-\mu(t)} d{t}}}\right).$ Arguing as in Case 1 of Theorem \ref{t3}, we conclude that $\lim_{n\rightarrow +\infty} P^{n}(w_{0})=0$.
		\end{proof}
		Combining Theorems \ref{t2}-\ref{t3} and Lemma \ref{lb1}, we can derive precise relationships between $g_{0}$ and the dynamical behavior of (\ref{eq1}). 
		\begin{itemize}
			\item If 
			\begin{equation*}
				\frac{\int_{0}^{T} \paren{a(t)-\mu(t)} d{t}}{\int_{0}^{\overline{T}}\xi(t)d{t}}<g_{0},
			\end{equation*}
			the origin is a global attractor, (see Lemma \ref{lb1}).
			\item If
			\begin{equation*}
				\frac{\int_{\overline{T}}^{T} a(t)d{t}-\int_{0}^{T}\mu(t) d{t}}{\int_{0}^{\overline{T}}\xi(t)d{t}}<g_{0}<\frac{\int_{0}^{T} \paren{a(t)-\mu(t)} d{t}}{\int_{0}^{\overline{T}}\xi(t)d{t}},
			\end{equation*}
			the origin is a local attractor (can be a global attractor) and there are, at most, two non-trivial $T$-periodic solutions, (see Theorem \ref{t3}).
			\item If
			\begin{equation*}
				g_{0}<\frac{\int_{\overline{T}}^{T} a(t)d{t}-\int_{0}^{T}\mu(t) d{t}}{\int_{0}^{\overline{T}}\xi(t)d{t}},
			\end{equation*}
			the origin is unstable and there is a non-trivial $T$-periodic solution that is a global attractor, (see Theorem \ref{t2}).
		\end{itemize}
		 Note that if one of the previous quantities appearing in the inequalities is negative, we obtain an impossible or redundant condition. For example, if $\int_{\overline{T}}^{T} a(t)d{t} < \int_{0}^{T}\mu(t) d{t}$, then the origin cannot be unstable for any $g_0 > 0$. Similarly, if $\int_{0}^{T} \paren{a(t)-\mu(t)} d{t} < 0$, the origin is a global attractor independently of the value of $g_0$. 
			
\begin{figure}
			\centering
			\subfigure[Backward bifurcation,  $P''(0) = 2.6682$ ]{
				\label{Figure2a}
				\includegraphics[width=0.47\textwidth]{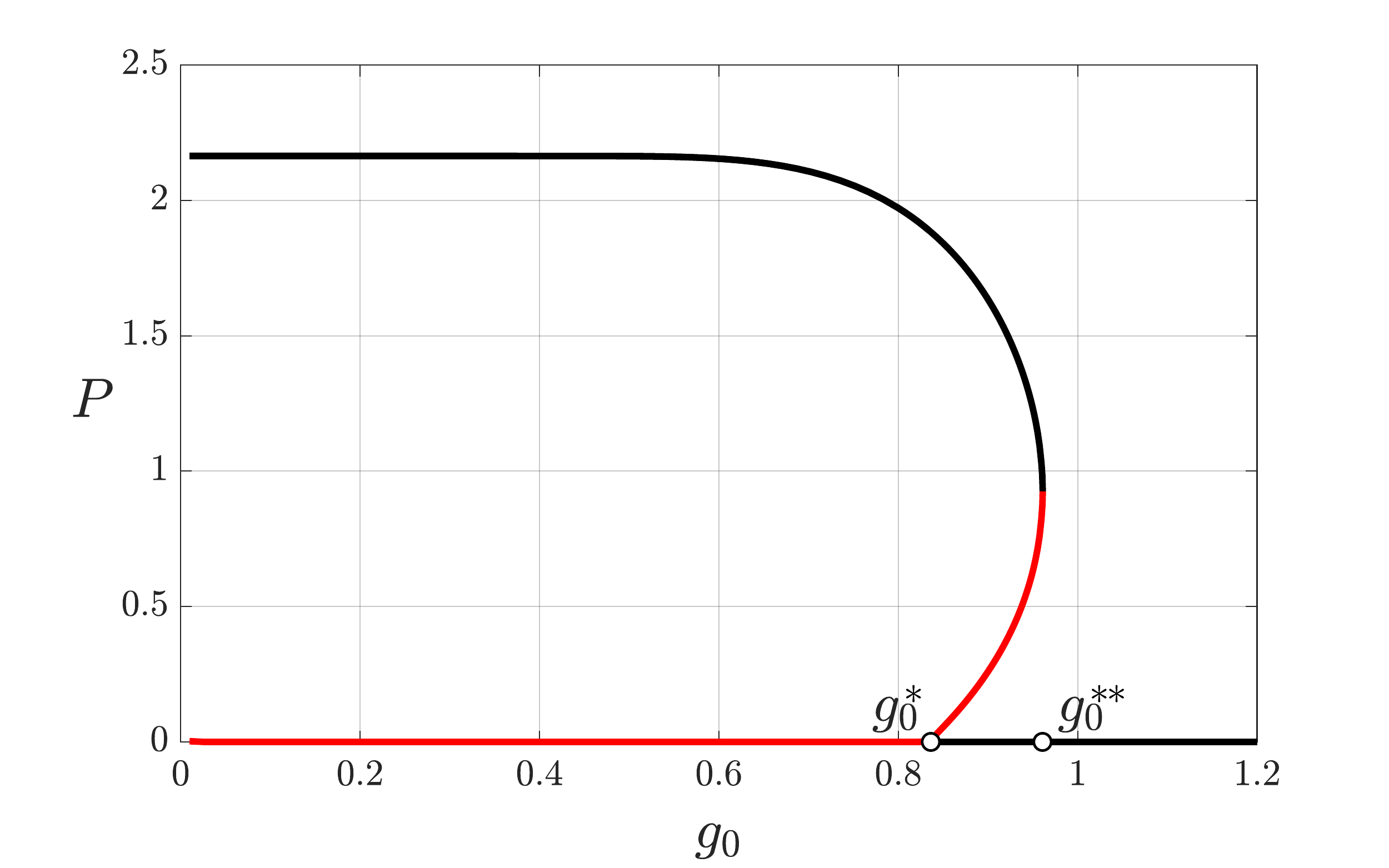}}
			\hfill
			\subfigure[Transcritical bifurcation in the limit case,  $P''(0) \simeq 0$]{
				\label{Figure2b}
				\includegraphics[width=0.47\textwidth]{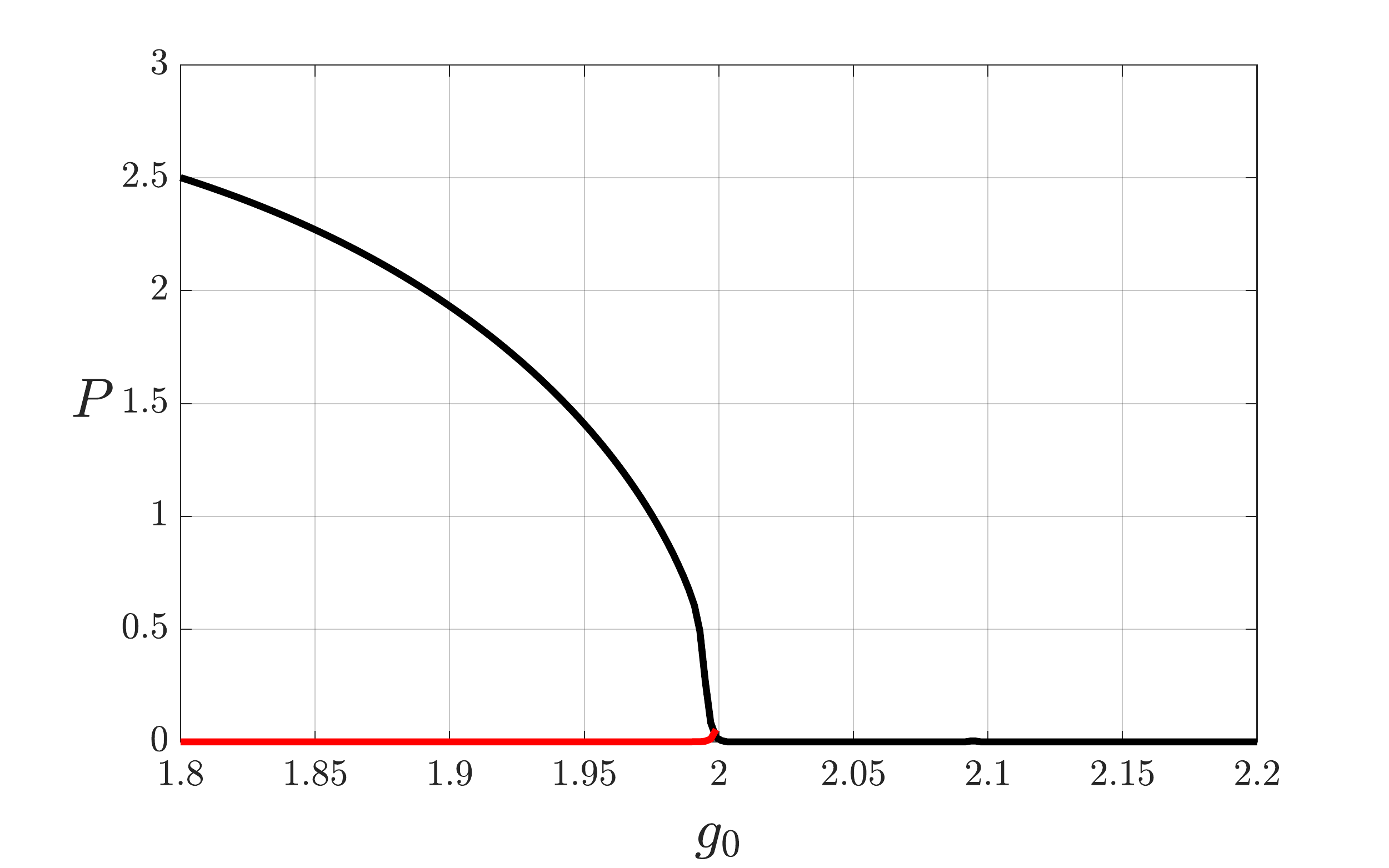}}
			\\
			\subfigure[Transcritical bifurcation, $P''(0) = -0.3056$]{
				\label{Figure2c}
				\includegraphics[width=0.45\textwidth]{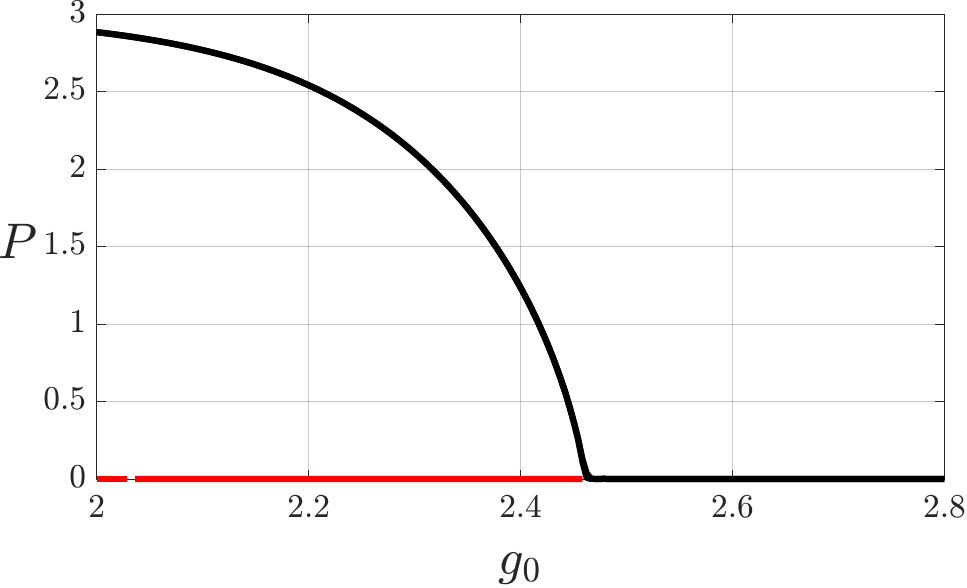}}
			
			\caption{ Bifurcation diagrams of the Poincaré map of (\ref{eq1}) with $g_{0}$ as bifurcation parameter. The black curve represents the stable fixed points of the Poincaré map, and the red one, the unstable fixed points. Parameters: (a) $a(t) = 4$ for $t \in [n,1/2)$ and $a(t) = 1$ for $t \in [n+1/2,n+1)$, with $n \in \mathbb{N} \cup \brac{0}$, $\xi(t) = 1$ for $t \in [n,1/2)$ and $\xi(t) = 0.2$ for $t \in [n+1/2,n+1)$, with $n \in \mathbb{N} \cup \brac{0}$, $\mu = 1$, $T = 14$, and $\overline{T} = 7$.  (b) $a = 4$, $\xi = 1$, $\mu = 1$, $T = 14$, and $\overline{T} = 7$.  (c) $a = 4$, $\xi = 1$, $\mu = 1$, $T = 14$, and $\overline{T} = 6.5$ }
			\label{fig: Bifurcaciones}
		\end{figure}

\begin{figure}
			\centering
			
			\subfigure[]{
				\includegraphics[width=0.45\linewidth]{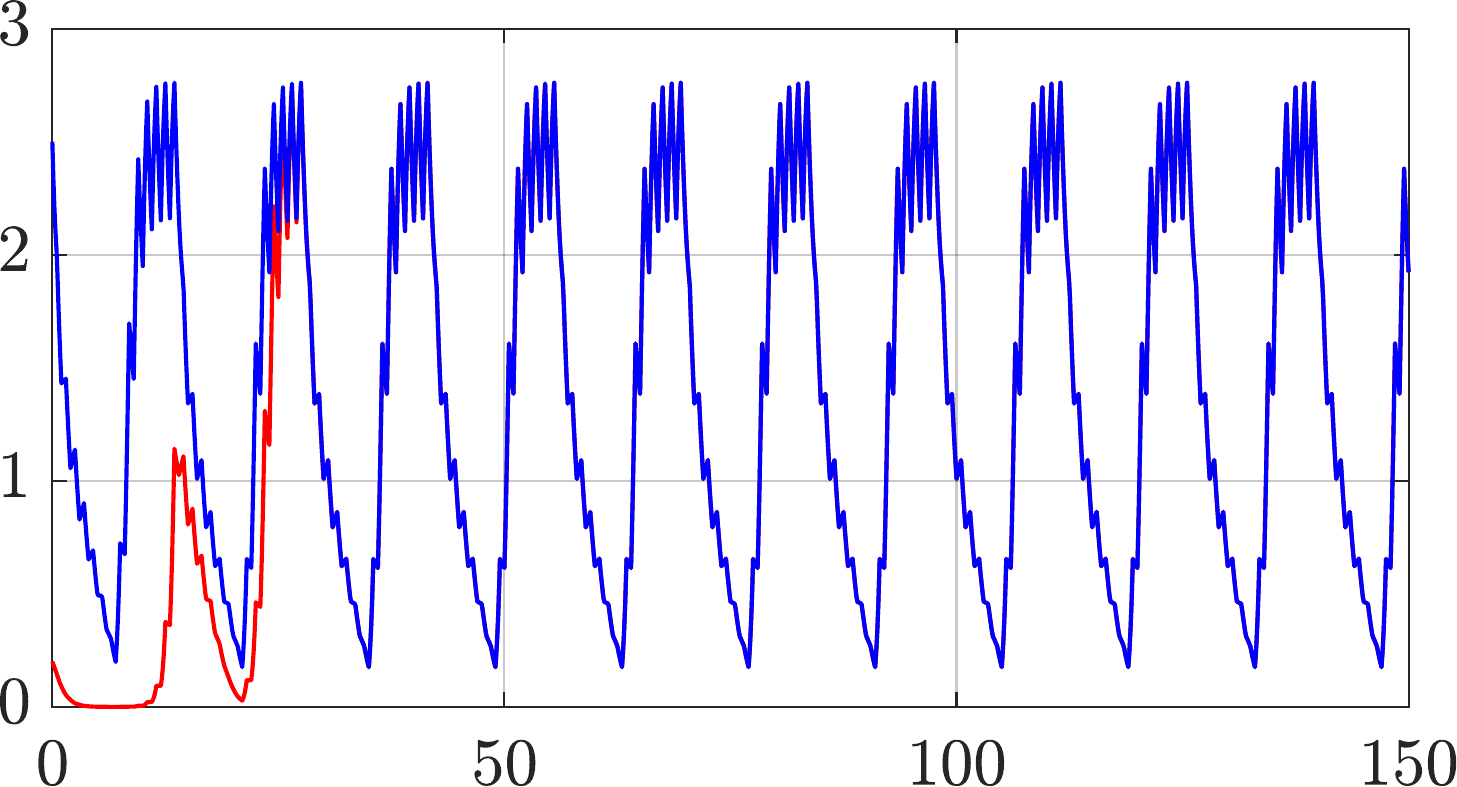}
			}
			\subfigure[]{
				\includegraphics[width=0.45\linewidth]{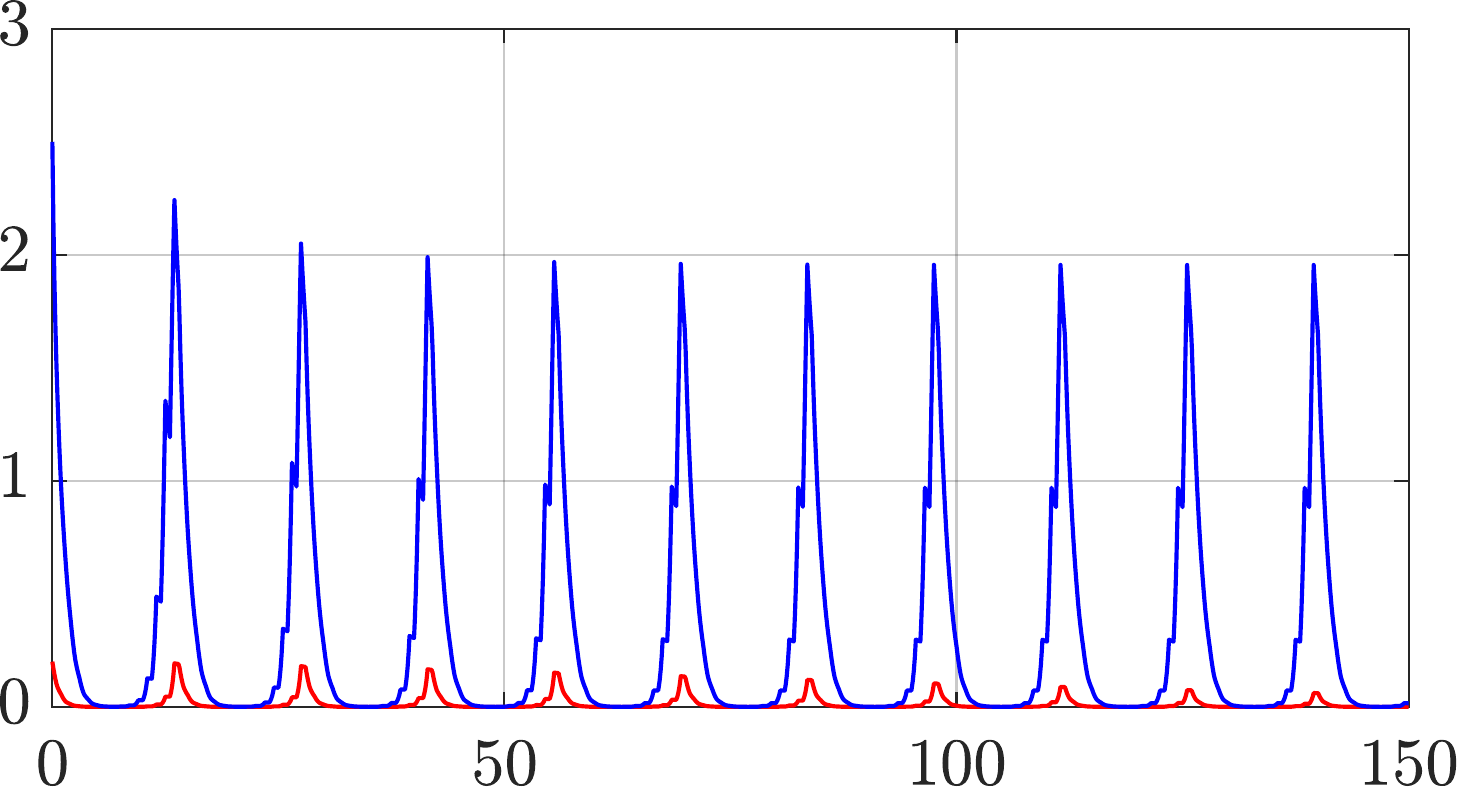}
			}
			\subfigure[]{
				\includegraphics[width=0.45\linewidth]{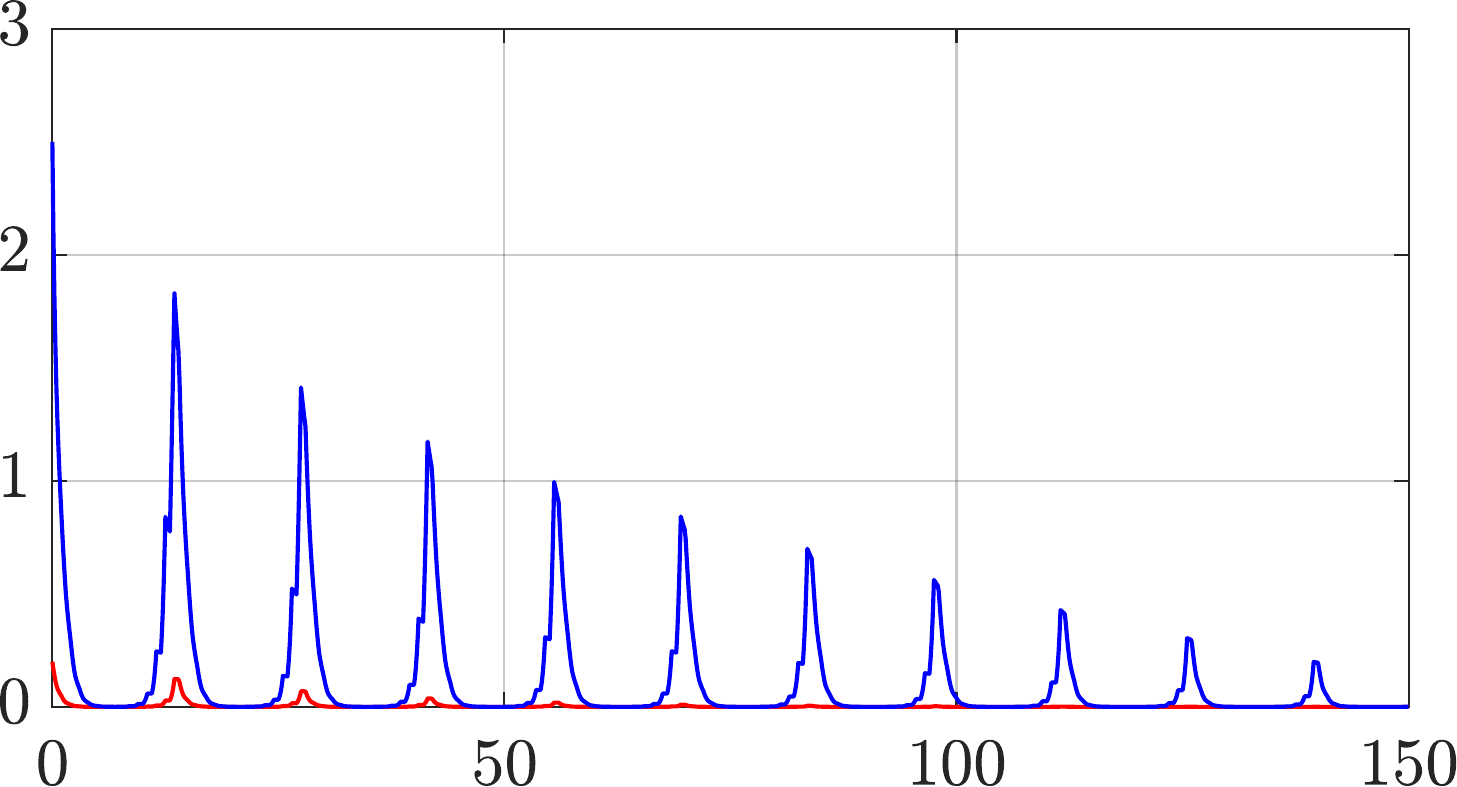}
			}
			\caption{Representation of the solutions of model (\ref{eq1}) for several values of $g_0$ and the initial condition. For the parameters employed in Figure \ref{Figure2a}, we illustrate the three biological scenarios discussed in the main text, taking initial conditions $0.2$ (red curves) and $2.5$ (blue curves). (a) $g_{0}=0.5$ (b) $g_{0}=0.9$ (c) $g_{0}=1$.}
			\label{fig:casosdiagramabifurcacion}
		\end{figure}

			Numerically analyzing the expected bifurcation diagrams of the Poincaré map of (\ref{eq1}) when $g_{0}$ is treated as a bifurcation parameter (Figure \ref{fig: Bifurcaciones}),  there are two dynamical scenarios: A usual backward bifurcation and a transcritical bifurcation. The value of $g_0$ where the origin turns from unstable to stable corresponds to $I_1 = 0$. Proposition \ref{pbif} classifies the bifurcations in terms of $P''(0)$. Specifically, if $P''(0)>0$, 
			a backward bifurcation emerges. On the other hand, if $P''(0)\leq 0$,  there is a transcritical bifurcation. In Figure \ref{fig:casosdiagramabifurcacion}, we draw several solutions to illustrate the convergence in the long term.
		

		\section{Extensions and further applications}
		In the derivation of model (\ref{intro2}), $T$ is the time between two consecutive interventions, e.g., a day or a month. On the other hand, the time dependence of the functions $a(t),\mu(t)$, and $\xi(t)$ represents the seasonal fluctuations of the environment. These functions are normally periodic but their period might differ from $T$ in many cases. Nevertheless, it seems natural to suppose that they are $\mathcal{T}$-periodic with $\mathcal{T}=n_{0}T$ for $n_{0}\in \mathbb{N}$. Under this condition, the analysis of (\ref{eq1}) is the same as that in Sections 2 and 3. There are mainly two differences: 1.-The Poincaré map is now defined as
		\begin{equation*}
			P:[0,+\infty)\longrightarrow P([0,+\infty))
		\end{equation*}
		\begin{equation*}
			P(w_{0})=w(\mathcal{T};w_{0}).
		\end{equation*}
		2.-(\ref{eq1}) has $n_{0}$ discontinuities in $(0,\mathcal{T})$ instead of one but the treatment is analogous, (see Theorem \ref{t1}).

		Model (\ref{intro2}) is a simple and versatile tool to predict the population dynamics of wild mosquitoes subject to the release of sterile mosquitoes. However, it could be oversimplifying in certain scenarios when other subtle reproduction and/or population growth aspects play a crucial role. To solve its possible limitations, several extensions of model (\ref{intro2}) have appeared in the literature. A key advantage is that our approach works for those extensions. Note that one of the main ingredients is that the sign of the third derivative of the inverse of the Poincaré map is strictly positive. As Theorem \ref{t1} indicated, this holds provided $\pdv[3]{F}{w} \paren{t,w} \leq 0$ for all $(t,w)\in [0,+\infty)^2$ and $\pdv[3]{F}{w} \paren{t,w} \not\equiv 0$ for some interval of $t$.
		
		Next, we discuss and analyze some variations of model (\ref{intro2}). As in Sections 2 and 3, we suppose that the birth $a(t)$, mortality $\mu(t)$, and competition rates $\xi(t)$ are $T$-periodic, of class $\mathcal{C}^{\infty}$, and strictly positive. In models 1 and 2, we consider the function $g$ defined in (\ref{g}). As mentioned above, the condition of the same period for all functions is not restrictive.\newline

		\noindent{\bf Model 1: A variation of the intraspecific competition.}\newline
		Wang {\it et al.} in \cite{Wang2024} analyzed the autonomous counterpart of the model
		\begin{equation}
			w' = \frac{a(t)w^2}{w + g(t)} \paren{1 - \eta \frac{w^2}{w + g(t)}} - \mu(t) w.
			\label{eq:Wang2024}
		\end{equation}
		The main difference regarding model (\ref{eq1}) is the term  $\paren{1 - \eta \frac{w^2(t)}{w(t) + g(t)}}$ which represents the intraspecific competition. Specifically, in (\ref{eq:Wang2024}), the expected number of newly born wild mosquito offspring at time $t$ is $\frac{a(t) w(t)}{w(t)+g(t)}$. Following \cite{Li2017}, the term $1 - \eta \frac{w^2}{w + g(t)}$ with $\eta>0$ represents the survival probability of newly born mosquitoes. Thus, the number of newly born mosquitoes surviving the intraspecific competition is 
		\begin{equation*}
			\frac{a(t)w^2}{w + g(t)} \paren{1 - \eta \frac{w^2}{w + g(t)}}.
		\end{equation*}
		See \cite{Wang2024} for a detailed discussion on the derivation of the model.  The third derivative of the system  is
		\begin{equation}
			\parcial[3]{F}{w}(t,w) = -\frac{6 a(t)g(t)^2 {\left(g(t)+w+4\eta(t) g(t)w\right)}}{{{\left(g(t)+w\right)}}^5 },
		\end{equation}
		so the extension of Theorem \ref{t1} is valid. \color{black}
		\newline
		Arguing as in Sections 2 and 3, we obtain the following result:
		\begin{theorem}
			Let $I_1 = \int_{\bar{T}}^{T} a(t) d{t} - \int_{0}^{T} \mu(t) d{t}$ and $I_2 = \int_{0}^{T} \paren{a(t) - \mu(t)} d{t}$.
			\begin{itemize}
				\item[i)] If $I_1 > 0$, the origin is unstable and there exists a non-trivial $T$-periodic solution $w(t;w_{0}^*)$ of model (\ref{eq:Wang2024}) that is an attractor in $(0,+\infty)$, i.e., for all $w_{0}\in (0,+\infty)$,
				\begin{equation*}
					\lim_{t\rightarrow +\infty}[w(t;w_{0})-w(t;w_{0}^*)]=0.
				\end{equation*}
				
				\item[ii)] If $I_1 < 0$ and $I_2 \leq 0$, \eqref{eq:Wang2024} does not have non-trivial $T$-periodic solutions and the origin is an attractor in $(0,+\infty)$, i.e., for each $w_{0}\in [0,+\infty)$, $\lim_{t\to+\infty}w(t;w_{0})=0$.
				
				\item[iii)] If $I_1 < 0$ and $I_2 > 0$,  the origin is a local attractor and there are, at most, two non-trivial $T$-periodic solutions of model \eqref{eq:Wang2024}. Moreover, if $a(t) > \mu(t)$ for all $t \in [0,\overline{T}]$,  $\lim_{t\rightarrow +\infty}w(t;w_{0})=0$ for all $w_{0}\in \left(0,\frac{-I_{1}g_{0}}{I_{2}e^{\int_{0}^{\overline{T}}\paren{a(t)-\mu(t)}d{t}}}\right)$.
			\end{itemize}
		\end{theorem}
		\noindent{\bf Model 2: Introducing imperfect Cytoplasmatic Incompatibility (CI).}\newline
		Cytoplasmatic Incompatibility (CI) is a biological mechanism of sterility of female mosquitoes associated with the release of Wolbachia-infected males. If CI occurs due to the mating of an uninfected female with an infected male, the female remains sterile for the rest of its life. Note that some of these matings may not lead to CI and the female reproduction follows a regular pattern. Introducing the possibility of not always having CI, that is, imperfect CI, Liu {\it et al.} in \cite{Liu2022} analyzed 
		\begin{equation}
			w' =w \left(\frac{a(t) w}{2(w + 2 g(t))} + \frac{a(t)(1 - s_h) g(t)}{w + 2g(t)} - \mu(t) - \xi(t)(w + g(t))\right)
			\label{eq:Liu2022}
		\end{equation}
		with $s_{h}\in [0,1]$ the intensity of CI. In (\ref{eq:Liu2022}), $\frac{a(t) w(t)}{2(w(t) + 2 g(t))}$ is the total number of newly born offspring at time $t$ produced by wild females mating with wild males whereas  $\frac{a(t)(1 - s_h) g(t)}{w(t) + 2g(t)}$ denotes the total number of newly born offspring at time $t$ produced by wild females mating with infected males. See \cite{Liu2022} for more details on the derivation of (\ref{eq:Liu2022}). In this case, the third derivative of the system is given by
		\begin{equation}
			\parcial[3]{F}{w}(t,w) = -\frac{12a(t)g(t)^2 s_h }{{{\left(2g(t)+w\right)}}^4 },
		\end{equation}
		hence once again the results in Theorem \ref{t1} can be extended to this model.\color{black}

		The dynamical picture of (\ref{eq:Liu2022}) is described in the next theorem.
		\begin{theorem}
			Let $I_1 = \int_{0}^{\overline{T}} \paren{\frac{a(t)}{2}(1-s_h) - \xi(t) g_0} d{t} + \int_{\overline{T}}^{T} \frac{a(t)}{2} d{t} - \int_{0}^{T} \mu(t) d{t}$ and \\$I_2 = \int_{0}^{\overline{T}} \paren{\frac{a(t)}{2}(1-s_h) - \xi(t) g_0} d{t} + \int_{0}^{T} \frac{a(t)}{2} d{t} - \int_{0}^{T} \mu(t) d{t}$.
			\begin{itemize}
				\item[i)] If $I_1 > 0$, the origin is unstable and there exists a non-trivial $T$-periodic solution $w(t;w_{0}^*)$ of model (\ref{eq:Liu2022}) that is an attractor in $(0,+\infty)$, i.e., for all $w_{0}\in (0,+\infty)$,
				\begin{equation*}
					\lim_{t\rightarrow +\infty}[w(t;w_{0})-w(t;w_{0}^*)]=0.
				\end{equation*}
				
				\item[ii)] If $I_1 < 0$ and $I_2 \leq 0$, \eqref{eq:Liu2022} does not have non-trivial $T$-periodic solutions and the origin is an attractor in $(0,+\infty)$, i.e., for each $w_{0}\in [0,+\infty)$, $\lim_{t\to+\infty}w(t;w_{0})=0$.
				
				\item[iii)] If $I_1 < 0$ and $I_2 > 0$,  the origin is a local attractor and there are, at most, two non-trivial $T$-periodic solutions of \eqref{eq:Liu2022}. Moreover, if $a(t)\paren{1 - \frac{s_h}{2}} - \mu(t) - \xi(t)g_0 > 0$ for all $t \in [0,\overline{T}]$, $\lim_{t\rightarrow +\infty}w(t;w_{0})=0$ for all $w_{0}\in \left(0,\frac{-2I_{1}g_{0}}{I_{2}e^{\int_{0}^{\overline{T}} \paren{a(t)\paren{1 - \frac{s_h}{2}} - \mu(t) - \xi(t)g_0} d{t}}}\right)$.
			\end{itemize}
		\end{theorem}
		
		\noindent {\bf Model 3: Using saturating release inputs.}\newline
		\noindent Zhang {\it et al.} in \cite{Zhang2023} modified model (\ref{intro2}) through a saturated release strategy. Inspired by the usual
		Holling type II functional response for predator-prey models, the input now reads as
		\begin{equation*}
			g(t) = \left\{\begin{array}{ll}
				\frac{b w(t)}{1+w(t)} & \mathrm{if } t\in {[i T,iT + \overline{T})},\\[4pt]
				0 & \mathrm{if } t\in {[iT + \overline{T},(i+1)T)}
			\end{array}\right.
		\end{equation*}
		with $b>0$ and $i \in \mathbb{Z}$. Moreover, they assumed that the released mosquitoes do not alter the intraspecific competition. Following these assumptions, we arrive at the switching model
		\begin{equation}\label{switch1}
			w' =\frac{a(t) w(t) (w(t) + 1)}{w(t) + 1 + b} - (\mu(t) + \xi(t) w(t))w(t) 
		\end{equation}
		for  $t \in {[i T,iT + \overline{T})}$,
		\begin{equation}\label{switch2}
			w' = a(t) w(t) - (\mu(t) + \xi(t) w(t))w(t) 
		\end{equation}
		for  $t\in {[iT + \overline{T},(i+1)T)}$ and $i \in \mathbb{Z}$. See \cite{Zhang2023} for more details on the derivation of the model. The third derivative is given by
		\begin{equation}
			\parcial[3]{F}{w}(t,w) = \left\{\begin{array}{ll}
				-\frac{6a(t)b{\left(b+1\right)}}{{{\left(b+w+1\right)}}^4 } & \mathrm{if } t\in {[i T,iT + \overline{T})},\\[4pt]
				0 & \mathrm{if } t\in {[iT + \overline{T},(i+1)T)} ,
			\end{array}\right.
		\end{equation}
		hence the results in Theorem \ref{t1} are applicable.\color{black}
		\newline
		Repeating the arguments in Sections 2 and 3, we characterize the dynamical behavior of (\ref{switch1})-(\ref{switch2}).
		
		\begin{theorem}
			Define $I_1 = \int_{\overline{T}}^{T} a(t) d{t} + \int_{0}^{\overline{T}} \frac{a(t)}{b+1} d{t} - \int_{0}^{T} \mu(t) d{t}$ and $I_2 = \int_{0}^{T} a(t) d{t} + \int_{0}^{\overline{T}} \frac{a(t)}{b+1} d{t} - \int_{0}^{T} \mu(t) d{t}$.
			\begin{itemize}
				\item[i)] If $I_1 > 0$, the origin is unstable and there exists a non-trivial $T$-periodic solution $w(t;w_{0}^*)$ of model (\ref{eq:Liu2022}) that is an attractor in $(0,+\infty)$, i.e., for all $w_{0}\in (0,+\infty)$,
				\begin{equation*}
					\lim_{t\rightarrow +\infty}[w(t;w_{0})-w(t;w_{0}^*)]=0.
				\end{equation*}
				
				\item[ii)] If $I_1 < 0$ and $I_2 \leq 0$, then \eqref{eq:Liu2022} does not have non-trivial $T$-periodic solutions and for each $w_{0}\in [0,+\infty)$, $\lim_{t\to+\infty}w(t;w_{0})=0$.
				
				\item[iii)] If $I_1 < 0$ and $I_2 > 0$, then the origin is a local attractor and there are, at most, two $T$-periodic solutions for \eqref{eq:Liu2022}. Moreover, if $a(t) + \frac{a(t)}{b+1} - \mu(t) > 0$ for all $t \in [0,\overline{T}]$ then $\lim_{t\rightarrow +\infty}w(t;w_{0})=0$ for all $w_{0}\in \left(0,\frac{-I_{1}(1+b)}{I_{2}e^{\int_{0}^{\overline{T}} \paren{a(t) + \frac{a(t)}{b+1} - \mu(t)} d{t}}}\right)$.
			\end{itemize}
		\end{theorem}
		
		\noindent {\bf Model 4: Introducing the Allee effect.}\newline
		It is broadly documented in theoretical ecology that many species exhibit the Allee effect, {\it i.e.}, a reduced per capita population growth rate at low densities, (see \cite{Schreiber2003} and the references therein).  If the Allee effect is mainly associated with the availability of mates, the birth rate is usually given by
		$\frac{a w}{w+1}$
		with $a>0$ the maximum growth rate. A usual model to describe the evolution of a population subject to the Allee effect is
		\begin{equation*}\label{inter1}
			w'=w\left(\frac{a(t)w}{w+1}-\mu(t)w-\xi(t)\right).
		\end{equation*}
		Chen {\it et al.} in \cite{Chen2023} studied the performance of  the release of Wolbachia-infected males through the switching model
		\begin{equation}\label{allee1}
			w'=\left(\frac{a(t)w}{(1+\alpha) w +1}-\mu(t)-\xi(t)w \right)
		\end{equation}
		for $t\in {[i T,iT + \overline{T})}$ and $i \in \mathbb{Z}$ with $\alpha>0$ and 
		\begin{equation*}\label{allee2}
			w'=\left(\frac{a(t)w}{ w +1}-\mu(t)-\xi(t)w \right)
		\end{equation*}
		for $t\in {[iT + \overline{T},(i+1)T)}$ and $i \in \mathbb{Z}$. In (\ref{allee1}), Wolbachia-infected males involve a reduction of the growth rate through the term $\alpha w$ in $\frac{a(t)w}{(1+\alpha) w +1}$.

		This model can be translated into model \eqref{eq1} after rewriting the parameters. Specifically,
		\begin{equation*}
			w'=w\left(\frac{\tilde{a}(t)w}{w+g(t)}-\xi(t)w-\mu(t)\right)
		\end{equation*}
		with
		\begin{equation*}
			\tilde{a}(t) =
			\left\{\begin{array}{ll}
				\frac{a(t)}{1 + \alpha} & \mathrm{if\ } t\in[i T,iT + \overline{T}),\\[4pt]
				a(t) & \mathrm{if\ } t\in[iT + \overline{T},(i+1)T) ,
			\end{array}\right.
			\qquad
			g(t) =
			\left\{\begin{array}{ll}
				\frac{1}{1 + \alpha} & \mathrm{if\ } t\in[i T,iT + \overline{T}),\\[4pt]
				1 & \mathrm{if\ } t\in[iT + \overline{T},(i+1)T) ,
			\end{array}\right.
		\end{equation*}
		for $i \in \mathbb{Z}$. Theorem \ref{t1} applies directly to the model. Since $g(t) \neq 0$ for every $t \geq 0$, \color{black} a key property of this model is that the origin is always a local attractor.
		
		\begin{theorem}
			Define $I = \int_{0}^{\overline{T}} \frac{a(t)}{\alpha + 1} d{t} + \int_{\overline{T}}^{T} a(t) d{t} - T$.
			\begin{itemize}	
				\item[i)] If $I \leq 0$, then \eqref{eq:Wang2024} does not have non-trivial $T$-periodic solutions and for each $w_{0}\in [0,+\infty)$, $\lim_{t\to+\infty}w(t;w_{0})=0$.
				
				\item[ii)] If $I > 0$, then the origin is a local attractor and there are, at most, two $T$-periodic solutions for \eqref{eq:Wang2024}. Moreover, if $a(t) > 1$ for all $t \in [0,\overline{T}]$ then $\lim_{t\rightarrow +\infty}w(t;w_{0})=0$ for all $w_{0}\in \left(0,\frac{-I_{1}}{I_{2}e^{\int_{0}^{\overline{T}} \paren{\frac{a(t)}{\alpha + 1} - 1} d{t}}}\right)$.
			\end{itemize}
		\end{theorem}

		\section{Discussion}
		Most organisms in nature are subject to some form of environmental seasonality, including temperature, photoperiod, precipitation, wind, and some human activities, (see \cite{El-morshedy2024,Lou2019,Ruiz-herrera2022,Gouagna2015,LeGoff2019,bellver2024dynamics} and the references therein).  It is clear from this list that seasonality is a notable aspect of the population dynamics of any species.  We stress that seasonality plays a key role in many invertebrates, such as ticks, mosquitoes, etc., due to diapause stages, (see \cite{El-morshedy2024,Lou2019,Ruiz-herrera2022,Hidalgo2018,Huestis2012,Denlinger2014} and the references therein). For these species,  overly cold winters or too dry seasons normally lead to reduced morphogenesis and physical activity. Consequently, neglecting seasonality might be unrealistic for the study of the dynamical behavior of mosquito populations.

		In this paper, mosquito demographic parameters (\textit{i.e.} birth rates and survival schedules) are time dependent to model the seasonal fluctuations of the environment.  Although theoretical ecologists broadly recognize its importance, many models frequently ignore seasonality. The main reason for this lack of interest is the mathematical complexity of the models. Generally speaking, non-autonomous models are much more difficult to study analytically. On the other hand, introducing seasonality can enlarge the possible dynamical behaviors of the models. This is particularly clear in classical epidemiological \cite{Barrientos2017,Ruiz-herrera2020,Ruiz-herrera2023} and predator-prey models \cite{Ruiz-herrera2012} where seasonal changes in the parameters alter the dynamical behavior from simple to chaotic dynamics.

		This paper aimed to compare and contrast the performance of the sterile insect technique when the demographic parameters of wild mosquitoes vary seasonally versus the autonomous case. To address this problem, we have proposed a methodology based on the sign of the third derivative of the inverse of the Poincaré map, together with subtle perturbative arguments. Our methodology was strongly motivated by some nice papers by Dueñas, Núñez, and Obaya \cite{duenas2023bifurcation,duenas2023critical,duenas2024critical,duenas2024generalized}. In comparison with those works, our paper has a double contribution: 1.- The geometric flavour of our analysis and the simplification of some results in the periodic case. 2.-The application in the context of the mosquito suppression models. Within this framework, we can compute various quantities that explicitly determine the dynamical behavior of the models. We emphasize that our arguments are completely different from those employed in  \cite{Wang2024, Yu2020a,Yu2020,Zheng2022,Zheng2023,Zheng2021b} and the references therein, where the demographic parameters are independent of time.

\begin{figure}
			\centering
			\subfigure[Solutions]{
				\label{Figure4a}
				\includegraphics[width=0.64\textwidth]{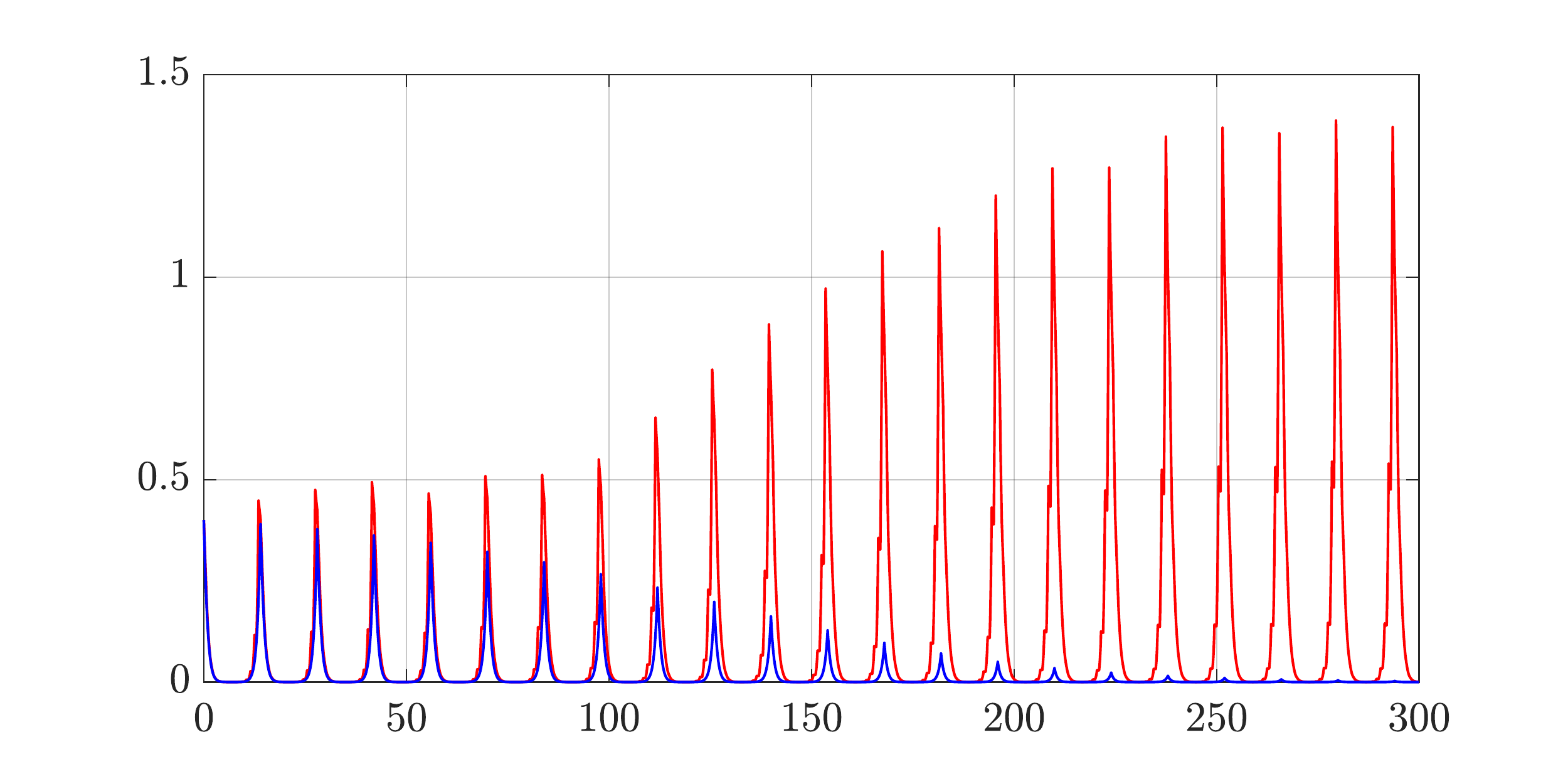}}
			\enspace
			\subfigure[Poincaré Map]{
				\label{Figure4b}
				\includegraphics[width=0.32\textwidth]{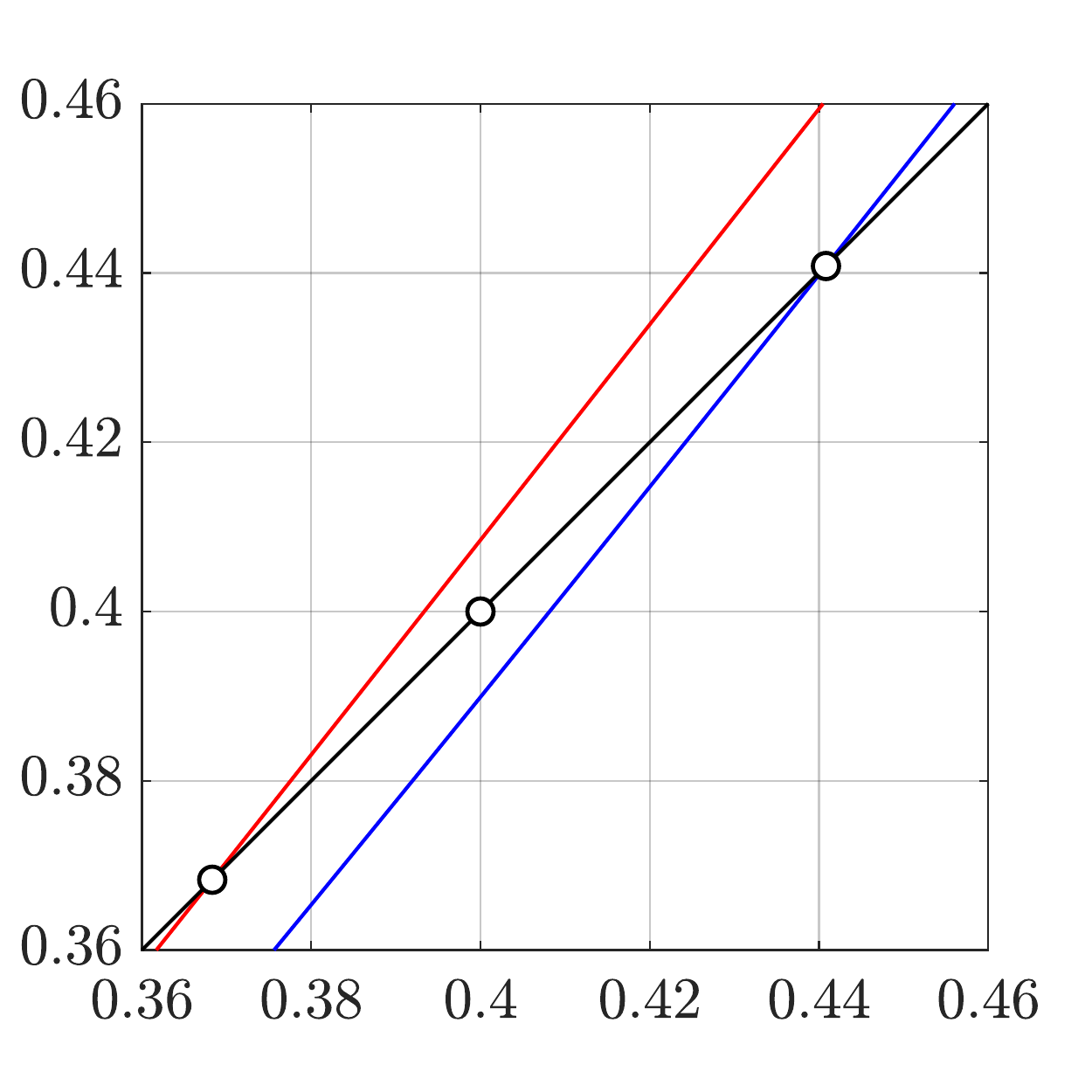}}
			
			\caption{ Comparison of the dynamical behavior of \eqref{eq1} with time-dependent parameters (red curves) and averaged parameters (blue curves). Fixed parameters $a(t) = 4$ for $t \in [n,1/2)$ and $a(t) = 1$ for $t \in [n+1/2,n+1)$, with $n \in \mathbb{N} \cup \brac{0}$, $\xi(t) = 1$ for $t \in [n,1/2)$ and $\xi(t) = 0.2$ for $t \in [n+1/2,n+1)$, with $n \in \mathbb{N} \cup \brac{0}$, $\mu = 1$, $T = 14$, $\overline{T} = 7$, and $g_0 = 0.75$. (a) Representation of the solution of both models with initial condition $w_{0}=0.4$. Observe that the solution associated with the model with averaged parameters (blue curve) tends to zero whereas the other one tends to a non-trivial $T$-periodic solution. (b) Representation of the Poincaré maps in the interval $(0.36,0.46)$. The unstable fixed points are different, so the basins of attraction of the origin are different. The highlighted points are the two unstable fixed points on the sides and the initial condition in the middle.}
			\label{Figure4}
		\end{figure}

\begin{figure}
			\centering
			\subfigure[Solutions]{
				\label{Figure5a}
				\includegraphics[width=0.64\textwidth]{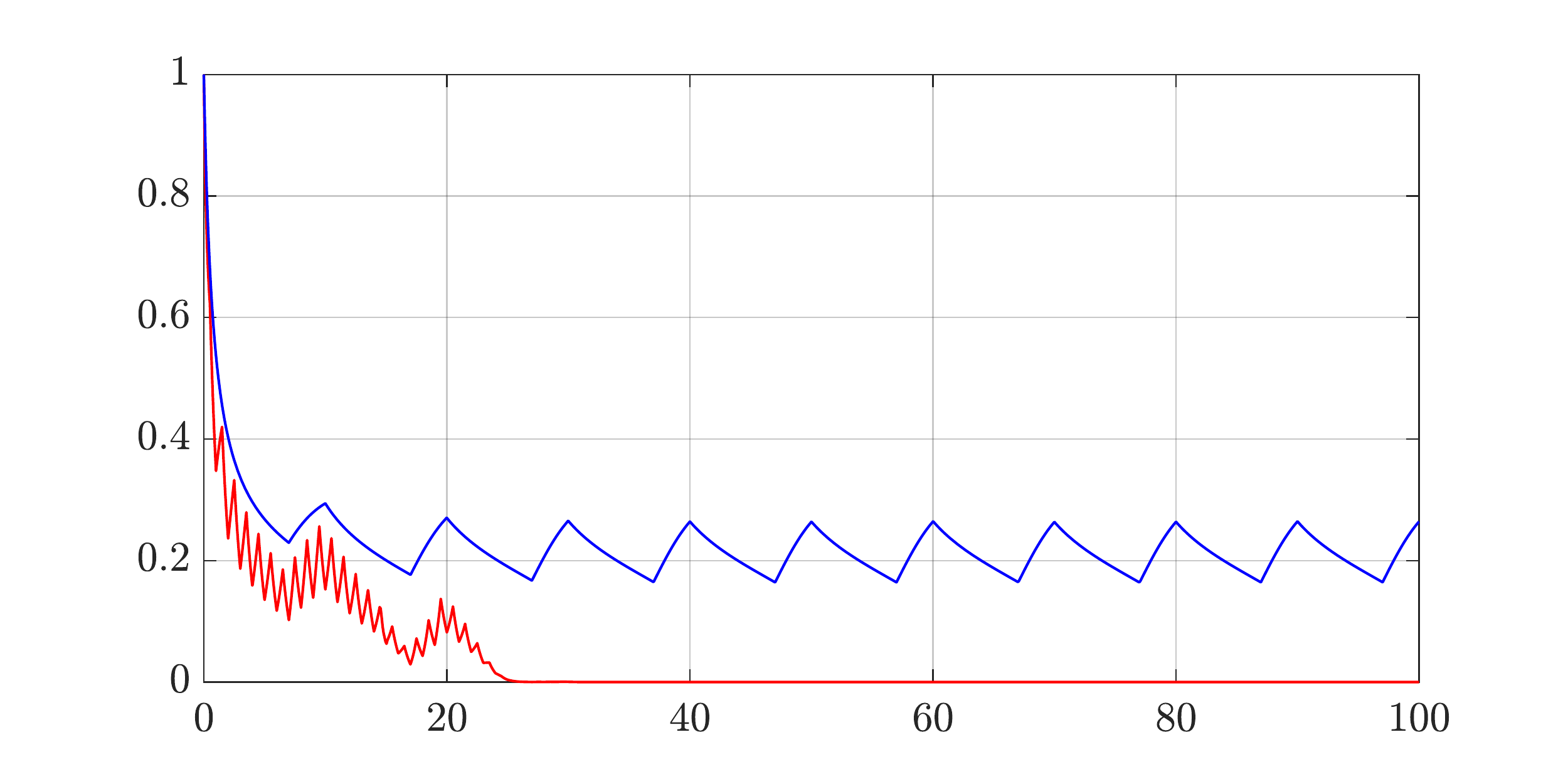}}
			\enspace
			\subfigure[Poincaré Map]{
				\label{Figure5b}
				\includegraphics[width=0.32\textwidth]{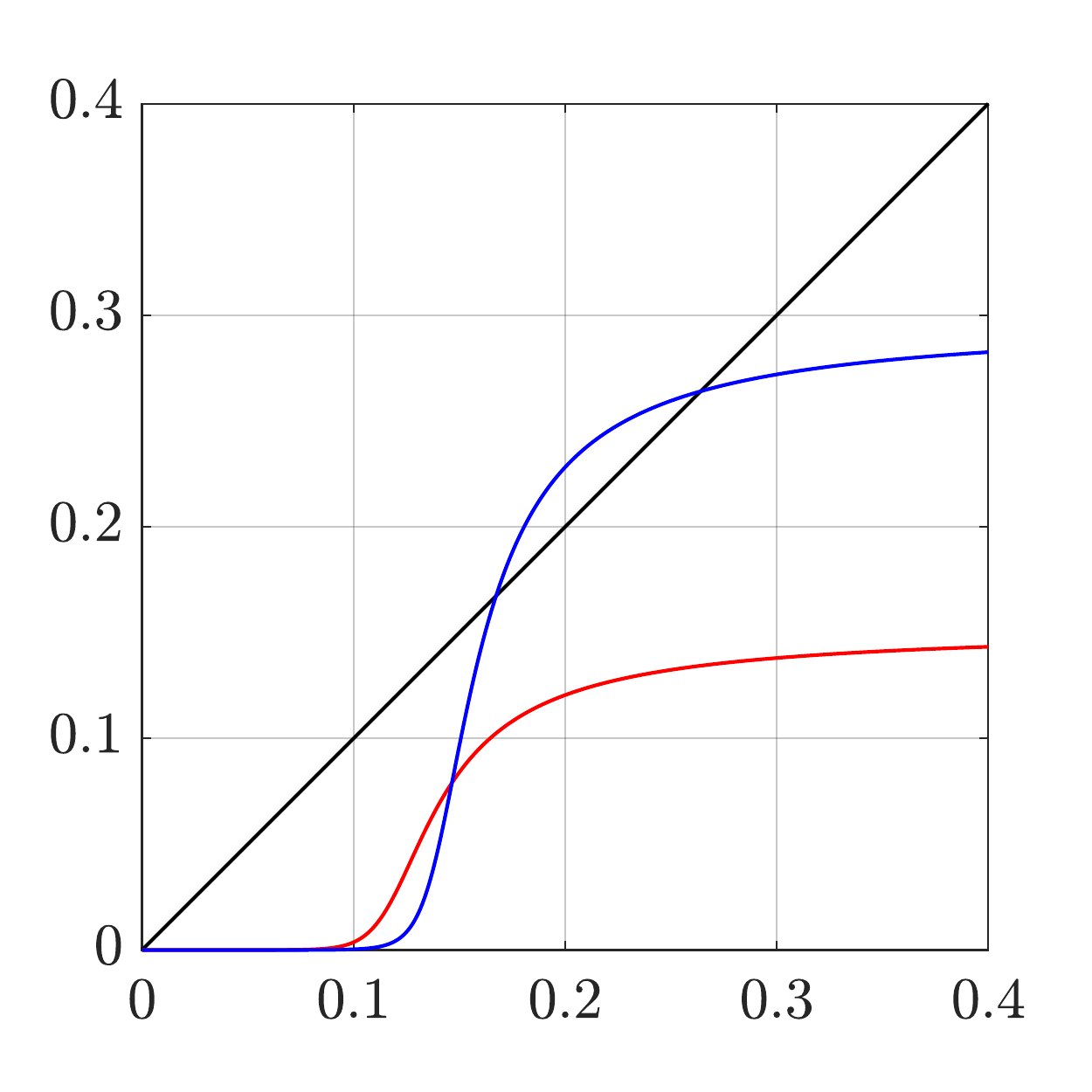}}
			
			\caption{ Comparison of the dynamical behavior of (\ref{eq:Wang2024}) with time-dependent parameters (red curves) and averaged parameters (blue curves). We have global attraction to the origin in the first case, and a bi-stability in the second one.
				Fixed parameters $a(t) = 4$ for $t \in [n,1/2)$ and $a(t) = 1$ for $t \in [n+1/2,n+1)$, with $n \in \mathbb{N} \cup \brac{0}$, $\eta(t) = a(t)/4$, $\mu = 2$, $T = 10$, $\overline{T} = 7$, and $g_0 = 0.03$. (a) Representation of the solution of both models with initial condition $1$.  (b) Representation of the Poincaré maps in the interval $(0,0.4)$.   }
			\label{Figure5}
		\end{figure}

		One conclusion of this work is that the stability of the trivial solution determines the dynamical behavior of the models. Roughly speaking, if the trivial solution is unstable, there exists a positive globally attractive periodic solution, while if the trivial solution is stable, there exists either global extinction or bi-stability with a non-trivial periodic solution. Therefore, seasonality does not create new dynamical behaviors in (\ref{intro2}). In other words, the results derived in \cite{Yu2020a,Yu2020} are valid when seasonality is introduced in the vital parameters. However, as Figures \ref{Figure4} and \ref{Figure5} indicated, there are strong differences between the autonomous and nonautonomous variants of the models.  Figures \ref{Figure4} and \ref{Figure5} compare \eqref{eq1} and (\ref{eq:Wang2024}) respectively with time dependent parameters and when considering the averaged parameters. Specifically, Figure \ref{Figure4} is an example in which both models exhibit a bi-stability. However, the basins of attraction of the origin are different.  
		On the other hand, Figure \ref{Figure5} is an example in which the dynamical behavior of both models is completely different. 
		
		\quad\newline\newline\noindent {\bf Acknowledgments}\newline
		We thank prof. R. Ortega for suggesting  the use of Lloyd's formulas. This work was supported by Ministerio de Ciencia, Innovación y Universidades through the grant n. PID2021-128418NA-I00.

		\appendix

		\section{Appendix}
		
		\subsection{Proof of Lemma \ref{l1}:}
		
		It is clear that $w'(0^+;w_0) < \widetilde{w}'(0^+;w_0)$ by \eqref{desi1}. Thus, there is $\delta > 0$ so that $w(t;w_0) < \widetilde{w}(t;w_0)$ for all $t \in (0,\delta)$. Assume, by contradiction, that there is a time $t_0 \in (0,\beta)$ so that
		\begin{equation}\label{eq:Lemma2.1-1}
			w(t_0;w_0) = \widetilde{w}(t_0;w_0)
		\end{equation}
		and
		\begin{equation}\label{eq:Lemma2.1-2}
			w(t;w_0) < \widetilde{w}(t;w_0) \quad \mathrm{for all } t \in (0,t_0) .
		\end{equation}
		Note that \eqref{eq:Lemma2.1-1} and \eqref{eq:Lemma2.1-2} imply that
		\begin{equation}\label{eq:Lemma2.1-3}
			w'(t_0^-;w_0) \geq \widetilde{w}'(t_0;w_0) .
		\end{equation}
		We distinguish between two cases:
		\begin{itemize}
			\item \textbf{Case 1}: $t_0 \notin \{i T+\overline{T},(i+1) T:i=0,1,2,\ldots\} $
			
			In this case, $w(t,w_0)$ is differentiable at $t_0$. Using that $w'(t_0;w_0) = F(t_0,w(t_0;w_0)) < G(t_0,\widetilde{w}(t_0;w_0)) = \widetilde{w}'(t_0;w_0)$, we obtain a contradiction with \eqref{eq:Lemma2.1-3}.
			
			\item \textbf{Case 2:} $t_0 \in \{i T+\overline{T},(i+1) T:i=0,1,2,\ldots\} $
			
			Define $\varepsilon = \widetilde{w}\paren{\frac{\delta}{2};w_0} - w\paren{\frac{\delta}{2};w_0}$. By continuous dependence on the initial conditions, we can take $w_0^* < w_0$ so that 
			\begin{equation}
				0 < \widetilde{w}(t;w_0) - \widetilde{w}(t;w_0^*) < \frac{\varepsilon}{2} \quad \mathrm{for all }t \in [0,t_0].
			\end{equation}
			The solutions $\widetilde{w}(t;w_0^*)$ and $w(t;w_0)$ satisfy that there are two instants $t_1,t_2$ with $0 < t_1 < t_2 < t_0$ such that 
			\begin{equation*}
				\widetilde{w}(t_1;w_0^*) = w(t_1;w_0) \quad \mathrm{and} \quad \widetilde{w}(t_2;w_0^*) = w(t_2;w_0),
			\end{equation*}
			and $\widetilde{w}(t;w_0^*) > w(t;w_0)$ for all $t \in (t_1,t_2)$. Repeating the arguments in Case 1 with initial time $t_1$ and initial condition $\widetilde{w}(t_1;w_0)$, we obtain a contradiction. We stress that taking $w_{0}^*$ close enough to $w_{0}$, we can guarantee that $t_{2}\not\in\{i T+\overline{T},(i+1) T:i=0,1,2,\ldots\}. $
		\end{itemize}
		
		\subsection{Computation of $P''(0)$}
		We recall that $P(w_{0})=w_{0}\Psi(w_{0})$ with 
		\begin{equation*}
			\Psi(w_{0})=\exp\left(\int_{0}^{\overline{T}} \paren{\frac{a(t) w(t;w_{0})}{w(t;w_{0})+g_{0}}-\mu(t)-\xi(t)(w(t;w_{0})+g_{0})} d{t}+\int_{\overline{T}}^{T} \paren{a(t)-\mu(t)-\xi(t)w(t;w_{0})} d{t}\right).
		\end{equation*}
		and $P$ is of class $\mathcal{C}^{\infty}$.
		Thus, $P''(0)=2\Psi'(0)$. Notice that 
		\begin{equation*}
			\Psi'(0)=\paren{\int_{0}^{\overline{T}} \paren{\frac{a(t)}{g_0} - \xi(t)} \parcial{w}{w_0}(t;0)  d{t} - \int_{\overline{T}}^{T} \xi(t) \parcial{w}{w_0}(t;0) d{t}} \Psi(0)
		\end{equation*}
		with
		\begin{equation*}
			\parcial{w}{w_0}(t;0)=e^{\int_{0}^{t}\frac{\partial}{\partial w}F(s,0)ds}
		\end{equation*}
		and 
		\begin{equation*}
			\frac{\partial F}{\partial w}(t,0) = \left\{\begin{array}{ll}
				- \mu(t) - \xi(t) g_0 & \mathrm{if } t\in  [i T,iT + \overline{T}),
				\\[10pt]
				a(t) - \mu(t) & \mathrm{if } t\in [iT + \overline{T},(i+1)T),
			\end{array}\right.
			\quad 
			i\in \mathbb{Z}.
		\end{equation*}
		Collecting this information, we arrive at

		\begin{align*}
			\Psi'(0) &= 
			\int_{0}^{\overline{T}} \paren{\frac{a(s)}{g_0} - \xi(s)} \exp\paren{-\int_{0}^{s} (g_0 \xi(t) + \mu(t)) d{t} } d{s} -
			\\[10pt]
			& - \int_{\overline{T}}^{T} \xi(s) \exp\paren{\int_{\overline{T}}^{s}a(t) d{t} - \int_{0}^{\overline{T}} g_0 \xi(t) d{t} - \int_{0}^{s} \mu(t) d{t} } d{s}.
		\end{align*}

		\color{black}
		\bibliographystyle{custombib}
		\bibliography{references}
		
	\end{document}